\documentclass[11pt]{article}

\usepackage{amsmath}
\usepackage{amssymb}
\usepackage{eucal}

\usepackage{stmaryrd}

\begin{document}

\baselineskip 16pt

\title{Finite groups in which generalized normality is a transitive  relation\thanks{Research was supported by  Ministry of
 Education of the Republic of Belarus (No. 20211328, No. 20211778). }}

\author{
Inna N. Safonova \\
{\small Department of Applied Mathematics and Computer Science,}\\
 {\small Belarusian State University,}\\
{\small Minsk 220030, Belarus}\\
{\small E-mail: in.safonova@mail.ru}\\
\\
Alexander N. Skiba\\
{\small Department of Mathematics and Technologies of Programming,}\\
{\small  Francisk Skorina Gomel State University,}\\
{\small Gomel 246019, Belarus}\\
{\small E-mail: alexander.skiba49@gmail.com}}

\date{}
\maketitle

\begin{abstract}   In this paper, we discuss some well-known results and some
 open problems of the theory  of $\sigma$-properties of a  group
 related to the study of generalized $T$-groups.

\end{abstract}


\let\thefootnote\thefootnoteorig

\section{Introduction}

Throughout this paper, all groups are finite and $G$ always denotes
a finite group; ${\cal L}(G)$ is the lattice of all subgroups of $G$.
 Moreover, 
 $\mathbb{P}$ is the set of all  primes and
 if
 $n$ is an integer, the symbol $\pi (n)$ denotes
 the set of all primes dividing $n$; as usual,  $\pi (G)=\pi (|G|)$, the set of all
  primes dividing the order of $G$;

A subgroup $A$ of $G$ is said to be \emph{quasinormal} (Ore)
 or \emph{permutable} (Stonehewer)
 in $G$
 if  $A$ permutes with every subgroup $H$ of $G$, that is, $AH=HA$; \emph{Sylow
 permutable} or \emph{$S$-permutable} \cite{prod, ill} if  $A$ permutes
 with all Sylow  subgroups of $G$.

A group  $G$ is said to be a \emph{$T$-group} if normality is a transitive
 relaion on $G$, that is,  if $H$ is a normal subgroup of $K$ and $K$ is a
 normal subgroup 
 of $G$, then $H$ is a normal subgroup of $G$. In other words, the group $G$
 is a $T$-group if every subnormal subgroup of $G$ is normal in $G$.

The   description of   $T$-groups  was first obtained by   Gasch\"{u}tz \cite{gasch} 
  for the soluble  case, and
  by Robinson in \cite{217}, for the general case.

The works \cite{gasch, 217} aroused great interest in the
 further study of $T$-groups
 and groups in which some conditions of generalized normality are  
 transitive ($PT$-groups, i.e. groups in which quasinormality
is transitive; $PST$-groups, i.e.
 groups, in which Sylow permutability is transitive, $P\sigma T$-groups, 
i.e.
 groups, in which $\sigma$-permutability is transitive (see Section 2 below), modularity
 is transitive   and    etc.).

However, there are many unsolved problems in this direction,
 and in this article we discuss some of them.

\section{$P\sigma T$-groups}

In what follows, $\sigma =\{\sigma_{i} \mid
 i\in I \}$ is some  partition of $\mathbb{P}$, that is, 
$\Bbb{P
}=\bigcup_{i\in I} \sigma_{i}$
 and $\sigma_{i}\cap
\sigma_{j}= \emptyset  $ for all $i\ne j$.

A \emph{$ \sigma $-property} of a group $G$ \cite{ProblemI, 1, 2, commun}   is any 
of  its property which does not depend on the  the choice of the partition  $ \sigma $
  of $ \mathbb {P}$. In other words, in the 
theory of $ \sigma $-properties of a group, we do not impose any restrictions on the 
partition $\sigma$ of $\mathbb{P}$.

Before continuing, recall some basic concepts of the theory  of $ \sigma $-properties 
of a group (see  \cite{ProblemI, 1, 2, commun}).

 If
 $n$ is an integer,
  $\sigma (n)= \{ \sigma_{i}\mid \sigma_{i} \cap \pi (n)\ne \emptyset\}$ 
 $\sigma (G)= \sigma (|G|)$.
 A group $G$ is said to be: \emph{$\sigma$-primary} if $G$ is a $\sigma _{i}$-group
 for some $i$;    \emph{$\sigma$-nilpotent} if $G$ is a direct product
 of $\sigma$-primary groups;  \emph{$\sigma$-soluble} if every chief factor of $G$ 
 is  $\sigma$-primary.

A subgroup $A$ of $G$ is said to
 be:

 (i) \emph{$\sigma$-subnormal}
  in $G$  if  there is a subgroup chain
 $$A=A_{0} \leq A_{1} \leq \cdots \leq
A_{n}=G$$  such that  either $A_{i-1}  
$ or
$A_{i}/(A_{i-1})_{A_{i}}$ is  ${\sigma}$-primary 
  for all $i=1, \ldots , n$;

 (ii) \emph{$\sigma$-seminormal
in $G$} (J.C. Beidleman) if $x\in N_{G}(A)$ for all $x\in G$ such that
 $\sigma (|x|)\cap \sigma (A)=\emptyset$;

 (iii) \emph{$\sigma$-permutable}
  in $G$  if either $A \trianglelefteq G$ or 
 $G$ is \emph{$\sigma$-full}, that is, $G$ has a Hall $\sigma _{i}$-subgroup 
for every $\sigma _{i}\in \sigma (G)$ and $A$ permutes with all such Hall
 subgroups of $G$.

{\bf Example  2.1.} (i) In the classical case
 $\sigma =\sigma ^{1}=\{\{2\}, \{3\}, \{5\}, \ldots
\}$    (we use here and below the notations in 
 \cite{alg12, 3})  a  subgroup $A$ of  $G$  is ${\sigma} ^{1}$-subnormal in $G$
 if and only if  $A$ is subnormal in $G$; $A$ is ${\sigma} ^{1}$-permutable  in $G$
if and only if  $A$ is Sylow permutable in $G$. A group $G$ is  ${\sigma} ^{1}$-soluble 
(${\sigma} ^{1}$-nilpotent) if and only if  $G$ is  soluble 
(respectively, nilpotent)

(ii) In the other classical case $\sigma =\sigma ^{\pi}=\{\pi,
\pi'\}$, $\pi\subseteq \Bbb{P}$,  a group $G$ is  ${\sigma} ^{\pi}$-soluble 
(${\sigma} ^{\pi}$-nilpotent) if and only if  $G$ is  $\pi$-separable 
(respectively,  $\pi$-decomposable, that is, $G=O_{\pi}(G)\times O_{\pi'}(G)$). 

A subgroup $A$ of  $G$  is
 ${\sigma} ^{\pi}$-subnormal  in $G$ if and only
 if  $G$ has   a subgroup chain
 $A=A_{0} \leq A_{1} \leq \cdots \leq
A_{n}=G$  such that  either $A_{i-1} \trianglelefteq A_{i}$, or
$A_{i}/(A_{i-1})_{A_{i}}$ is  ${\pi}$-group, or
 $A_{i}/(A_{i-1})_{A_{i}}$ is  ${\pi}'$-group  for all $i=1, \ldots , n$.

In this case we  say that $A$ is
 \emph{$\pi, \pi'$-subnormal} in $G$. 

A subgroup $A$ of  $G$  is   ${\sigma} ^{\pi}$-permutable   in $G$ if and only
 if  $G$ has a  Hall $\pi$-subgroup and a Hall $\pi'$-subgroup  and $A$ permutes with all 
such Hall subgroups of $G$.
In this case we  say that $A$ is
 \emph{$\pi, \pi'$-permutable} in $G$.

(iii) In fact, in the theory of $\pi$-soluble groups ($\pi= \{p_{1}, \ldots , p_{n}\}$)
 we deal with the  partition
$\sigma =\sigma ^{1\pi }=\{\{p_{1}\}, \ldots , \{p_{n}\}, \pi'\}$ of $\Bbb{P}$.  
A group $G$ is  ${\sigma} ^{1\pi}$-soluble 
(${\sigma} ^{1 \pi}$-nilpotent) if and only if  $G$ is  $\pi$-soluble 
(respectively,  $\pi$-special \cite{alg12, 3}, that is,
 $G=O_{p_{1}}(G)  \times   \cdots  \times O_{p_{n}}(G) \times O_{\pi'}(G)$).

 A subgroup $A$ of
 $G$ is $\sigma^{1\pi }$-subnormal in $G$ if and only if  $G$ has 
 a subgroup chain
 $A=A_{0} \leq A_{1} \leq \cdots \leq
A_{n}=G$  such that  either $A_{i-1} \trianglelefteq A_{i}$  or
$A_{i}/(A_{i-1})_{A_{i}}$ is  ${\pi}'$-group   for all $i=1, \ldots , n$.

In this case we say that $A$ is \emph{$\pi$-subnormal} in $G$.

In fact, the appearance of the theory of  $ \sigma $-properties 
of a group was
 connected chiefly with attempts to solve the following very dificult problem.

{\bf Question 2.2} (See Question in \cite{1}).
 {\sl What is the structure of  a $\sigma$-full group $G$ in which  
$\sigma$-permutability is transitive on $G$, that is, 
 if $H$ is a $\sigma$-permutable subgroup of $K$ and $K$ is a $\sigma$-permutable
 subgroup 
 of $G$, then $H$ is a  $\sigma$-permutable subgroup of $G$?
}

This problem turned out to be difficult even in the $\sigma$-soluble 
 case, its
 solution required the development of many aspects of the
 theory of $\sigma$-properties of a
 group. The theory of $\sigma$-soluble $P\sigma T$-groups
 was mainly developed in the papers 
 \cite{1, alg12, adar, 4, 2021, 444, 555, Ad} and the
following theorem (which, in fact, is the main result of the papers \cite{1, alg12})
 is the key result in this direction.

{\bf Theorem 2.3 } (See  Theorem A in
 \cite{alg12}).   {\sl If  $G$ is a  $\sigma$-soluble
 $P\sigma
T$-group and  $D=G^{\frak{N_{\sigma}}}$,
 then
    the following conditions hold:}

(i) {\sl $G=D\rtimes M$, where $D$   is an abelian  Hall
 subgroup of $G$ of odd order, $M$ is $\sigma$-nilpotent  and  every element of $G$ induces a
 power automorphism in $D$;  }

(ii) {\sl  $O_{\sigma _{i}}(D)$ has
a normal complement in a Hall $\sigma _{i}$-subgroup of $G$ for all $i$.}

{\sl Conversely, if  Conditions (i) and (ii) hold for  some subgroups $D$ and $M$ of
 $G$, then $G$ is  a  $\sigma$-soluble
 $P\sigma T$-group.}

In this theorem,  $G^{\frak{N_{\sigma}}}$  is the \emph{$\sigma$-nilpotent
 residual} of $G$,
 that is,  the intersection of all normal subgroups $N$ of $G$ with 
$\sigma$-nilpotent quotient $G/N$.

In the case
 $\sigma =\sigma ^{1}=\{\{2\}, \{3\}, \{5\}, \ldots
\}$, we get from Theorem 2.3 the following known result.

{\bf Corollary 2.4} (Agrawal  \cite[Theorem 2.3]{Agr}). {\sl Let 
$D=G^{\frak{N}}$  be 
 the nilpotent residual of $G$.  If  $G$ is a soluble $PST$-group,
 then  $D$  is  an
abelian  Hall subgroup of $G$ of odd  order and every element
 of $G$ induces a power automorphism 
in  $D$.  }

In order to consider some further
 applications of this Theorem 2.3, we  introduce the following concepts.

{\bf Definition 2.5.} Let $\mathfrak{X}$ be a class of groups. Suppose that
 with  each group $G\in \mathfrak{X}$ we
 associate some system of 
  its subgroups $\tau(G)$. Then we say that $\tau$ is a \emph{subgroup
 functor}  in the sense of Skiba 
 \cite{I} or $\tau$ is a \emph{Skiba subgroup functor} on $\mathfrak{X}$ \cite{II}
 if the following conditions are met:

 (1) $G\in\tau(G)$ for any group $G$;

(2) for any epimorphism $\varphi :A\mapsto B$, where
$A, \in \mathfrak{X}$,    and for any groups $H\in\tau(A)$ and $T\in\tau(B)$,
 we have $H^{\varphi}\in \tau(B)$ and $T^{{\varphi}^{-1}}\in \tau(A)$.

In what follows,  $\tau$ is some  subgroup
 functor on $\mathfrak{X}$  in the sense of Skiba.

If  $A\in \tau(G)$, then we say that $A$ is a \emph{$\tau$-subgroup} of $G$.
If  $\mathfrak{X}$ is the class of all groups, then instead of "subgroup
 functor on $\mathfrak{X}$" we will simply say "subgroup
 functor".

 We say also that a subgroup functor $\tau$ on $\mathfrak{X}$  
 is \emph{$\sigma$-special} 
 if  for  any  group $G\in \mathfrak{X}$  
the following three conditions  hild:

(*)  Each of  $\sigma$-subnormal $\tau$-subgroups of $G$ is
 $\sigma$-permutable in $G$, and

(**) $\langle A, B \rangle\in \tau (G)$ for any two $\sigma$-subnormal
  subgroups  $A, B \in \tau (G)$ of $G$,

(***)  If  $G=D\rtimes M$ is a $\sigma$-soluble $P\sigma T$-group, where 
 $D=G^{\frak{N_{\sigma}}}$ and  $A$ 
is a $\sigma$-primary $\sigma$-subnormal subgroup of $G$ such that  
  $A \in \tau (M)$, then  $A \in \tau (G)$.

{\bf Lemma 2.6.}  {\sl Suppose that $G=D\rtimes M$ is a
 $\sigma$-soluble $P\sigma T$-group, where 
 $D=G^{\frak{N_{\sigma}}}$. If  $A$ 
is a $\sigma$-primary $\sigma$-subnormal subgroup of $G$ such that  
  $A \leq M$, then $D\leq C_{G}(A)$. }

{\bf Proof.}  Let $A$ be a $\sigma _{i}$-group  and $x$  an element 
 in $D$ of prime power  order $p^{n}$.  The group $G$ is $\sigma$-soluble,
 so $G$ has a Hall $\sigma _{k}$-subgroup $H_{k}$ for all $k$  
 by Theorem B in \cite{2}.   In view of Theorem
 2.3, $H_{k}=O_{\sigma _{k}}(D)\times S_{i}$.

  Since  $A$ is 
 $\sigma$-subnormal in $G$, $A \leq H_{i}$ by Lemma 2.30(7) below. On the other
 hand, since 
$A \leq M$, $A\cap D=1$. Therefore  $A=(A\cap O_{\sigma _{i}}(D)) \times (A\cap S_{i})$
 since 
$O_{\sigma _{i}}(D)$ is a Hall subgroup of $D$ and of $G$, 
 so  $A\leq S_{i}$ and hence $ O_{\sigma _{i}}(D)\leq C_{G}(A)$. Now, let $k\ne i$.
 Then
$A$ is a Hall $\sigma _{i}$-subgroup of $V:=O_{\sigma _{k}}(D)A$ and $A$ is
 $\sigma$-subnormal
 in $V$ by Lemma 2.30(1) below, so $V=O_{\sigma _{k}}(D)\times A$  by 
Lemma 2.30(6) below
 and hence 
 $ O_{\sigma _{k}}(D)\leq C_{G}(A)$.  

The lemma is proved.

{\bf Lemma 2.7.}  {\sl   Let   $N\leq A$ be subgroups of  a $\sigma$-full
 group $G$, where    $N$ is normal in $G$.  Suppose that $ \{\sigma _{1}, \ldots ,
 \sigma _{t}\}=\sigma (G)$ 
  and $H_{i}$ is a Hall  $\sigma _{i}$-subgroup of $G$   for 
all $i=1, \ldots, t$.}

(1) {\sl If   $AH_{i}^{x}=H_{i}^{x}A$ for all $i=1, \ldots, t$ and all $x\in G$, then $A$ 
is $\sigma$-permutable in $G$.}
 (see Proposition 1.1 in \cite{2019}).

(2) {\sl   $A/N$ is  $\sigma$-permutable in $G/N$ if and only if  $A$ is 
$\sigma$-permutable in $G$}.

{\bf Proof.}  First note that 
$H_{i}N/N $ is a Hall $\sigma _{i}$-subgroup of $G/N$ for
 every $\sigma _{i}\in \sigma (G/N)$, so  $G/N$ is  $\sigma$-full. 

If   $AH_{i}^{x}=H_{i}^{x}A$ for all $i=1, \ldots, t$ and all $x\in G$, then 
$$AH_{i}^{x}/N=  (A/N)(H_{i}N/N)^{xN}=(H_{i}N/N)^{xN}(A/N)=H_{i}^{x}A/N$$
 for all $i=1, \ldots, t$ and all $xN\in G/N$. Therefore if 
$A$ is 
$\sigma$-permutable in $G$, then $A/N$ is  $\sigma$-permutable in $G/N$ by Part (1). 

Similarly, if   $A/N$ is  $\sigma$-permutable in $G/N$, then  $A$ is 
$\sigma$-permutable in $G$  by Part (1).

The lemma is proved.

{\bf Example 2.8.  }    Let $\mathfrak{X}$ be the  class of of all  $\sigma$-full
 groups.

 (1)  Let, for any $\sigma$-full group $G$,  $ \tau (G)$ be the set of all 
$\sigma$-permutable  subgroups of $G$.  Then, in view of Lemma 2.7(2), 
$\tau$ is a subgroup
 functor in the sense of Skiba and, clearly  Condition (*) holds  for  any 
 group $G\in \mathfrak{X}$. Moreover, in view of  \cite[A, 1.6(a)]{DH}, 
Condition  (**)
holds for any  group $G$.  Finally, 
 Conditions (***) holds in any $P\sigma T$-group by Theorem 2.3.

(2)  Recall that a  subgroup $M$ of $G$  is said to be  \emph{ modular} in $G$
   if $M$ is a modular element  (in the sense of
 Kurosh \cite[p. 43]{Schm})  of the  lattice ${\cal L}(G)$,  that is,

(i) $\langle X,M \cap Z \rangle=\langle X, M \rangle \cap Z$ for all $X \leq G, Z \leq
 G$ such that $X \leq Z$, and

(ii) $\langle M, Y \cap Z \rangle=\langle M, Y \rangle \cap Z$ for all $Y \leq G, Z \leq
 G$ such that  $M \leq Z$.

Let, for any $\sigma$-full group $G$,  $ \tau (G)$ be the set of all 
modular  subgroups of $G$.  Then, in view of  
\cite[p. 201, Properties (3), (4)]{Schm}, 
$\tau$ is a subgroup
 functor in the sense of Skiba.

Now, we show that  the functor $\tau$ is   $\sigma$-special. Indeed, if $A$ is a
 $\sigma$-subnormal modular subgroup of  a $\sigma$-full group $G$, then 
  $A$ is  $\sigma$-permutable in $G$ by Theorem 3.3 (i) below, so Condition (*)
 holds  for $G$.  Next let $A$ and $B$ be  $\sigma$-subnormal modular subgroups of $G$.
Then 
$\langle A, B \rangle$ is   modular in $G$ by \cite[p. 201, Property (5)]{Schm}, so
 Condition (**)   holds  for $G$.

Finally, suppose that  $G=D\rtimes M$ is a $\sigma$-soluble $P\sigma T$-group, where 
 $D=G^{\frak{N_{\sigma}}}$, and  let $A$ 
be a $\sigma$-primary $\sigma$-subnormal subgroup of $G$ such that  
  $A \in \tau (M)$, that is, $A$ is modular in $M$.  We show
 that in this case
 we have  $A \in \tau (G)$, that is, $A$ is modular in $G$.

 In view of Lemma 5.1.13 in \cite{Schm}, it is enough to show that if $x$
 is an element of $G$ of prime
 power order  $p^{n}$, then $A$ is modular in $\langle x, A \rangle$.
 
If $x\in D$, it is clear. Now assume that $x\not \in D$ and so $x\in M^{d}$ for
 some $d\in D$  since $M$ is a Hall subgroup of $G$. But $A$ is modular in $M$ and 
so $A$ is modular in $M^{d}$ since $A^{d}=A$ by Lemma 2.6. Therefore 
$A$ is modular in $\langle x, A \rangle$.   
Hence   Condition (***) holds for $G$,  so $\tau$ is a 
$\sigma$-special subgroup functor on $\mathfrak{X}$ .

(3)  Let, for any group $G$,  $ \tau (G)$ be the set of all 
normal subgroups of $G$.  Then   a subgroup functor $\tau$ is 
  $\sigma$-special (see Part (2)).

{\bf Lemma 2.9 } (See  Corollary 2.4 and Lemma 2.5  in \cite{1}).  {\sl
The class   of all  $\sigma$-nilpotent groups
 ${\mathfrak{N}}_{\sigma}$  is closed under taking
products of normal subgroups, homomorphic images and  subgroups. 
 Moreover, if $E$ is a normal subgroup of $G$ and  $E/(E\cap \Phi (G))$ 
  is $\sigma$-nilpotent, then $E$ 
is $\sigma$-nilpotent. }

{\bf Lemma 2.10 } (See  Proposition 2.3  in \cite{1}).  {\sl
A group $G$ is   $\sigma$-nilpotent if and only if every subgroup of $G$ is 
$\sigma$-subnormal in $G$.
  }

Now we prove the following result.

{\bf Theorem 2.11. }  {\sl  Suppose that $G$ is a $\sigma$-soluble group with
 $D=G^{\frak{N_{\sigma}}}$ and let 
 $\tau$ be  a  $\sigma$-special subgroup fuctor on the set of all $\sigma$-full groups
 $\mathfrak{X}$. 
If every $\sigma$-subnormal subgroup of $G$ is  $\tau$-subgroup of $G$, then $G$ is a   
  $P\sigma
T$-group and  the following conditions hold:}

(i) {\sl $G=D\rtimes M$, where $D$   is an abelian  Hall
 subgroup of $G$ of odd order, $M$ is a $\sigma$-nilpotent  group with $U\in \tau (M)$
 for all subgroups $U$ of $M$, and  every element of $G$ induces a
 power automorphism in $D$;  }

(ii) {\sl  $O_{\sigma _{i}}(D)$ has
a normal complement in a Hall $\sigma _{i}$-subgroup of $G$ for all $i$.}

{\sl Conversely, if  Conditions (i) and (ii) hold for  some subgroups $D$ and $M$ of
 $G$, then  every $\sigma$-subnormal subgroup of $G$ belongs to $\tau (G)$.}

{\bf Proof.}    First assume, arguing 
by contradiction, that Conditions
 (i) and   (ii) hold for some subgroups $D$ and $M$ of $G$ but $G$ has a 
 $\sigma$-subnormal
subgroup $U$   such that $U\not \in \tau (G)$. Moreover,
  we can assume that  $G$ is a counterexample with $|G|+|U|$ minimal.  Then  
 $U_{0} \in \tau (G)$ 
for every 
$\sigma$-subnormal  subgroup $U_{0}$ of $G$ such that  $U_{0} < U$.

(1)  $U\cap D=1$.   

Indeed, assume that $V:=U\cap D \ne 1$.  Then $G/V=(D/V)\rtimes (MV/V)$, where 
  $$D/V=G^{\frak{N_{\sigma}}}/V=
(G/V)^{\frak{N_{\sigma}}}$$    is an abelian  Hall
 subgroup of $G/V$ of odd order, $MV/V\simeq M$ is a $\sigma$-nilpotent
  group in which every subgroup
belongs to  $\tau (MV/V)$ and  every element of $G/V$ induces a
 power automorphism in $D/V$.  It is clear also that 
$O_{\sigma _{i}}(D/V)$ has
a normal complement in a Hall $\sigma _{i}$-subgroup of $G/V$ for all $i$. Therefore 
Conditions
 (i) and   (ii) hold for the  subgroups $D/V$ and $MV/V$ of $G/V$. Therefore 
$U/V \in \tau (G/V)$ 
for the  
$\sigma$-subnormal  subgroup $U/V$ of $G/V$ by the choice of $G$, so $U \in \tau (G)$ by
 the definition of the subgroup functor  $\tau$.  This contradiction shows that we have 
(1). 
 
(2) {\sl $U$ is a $\sigma _{i}$-group for some $i$.}

In view of Lemma 2.9 and Claim (1), $U\simeq UD/D$ is $\sigma$-nilpotent.
 Then $U$ is the 
direct product of some $\sigma$-primary  non-identity
 groups $U_{1}, \ldots , U_{t}$. If $U_{1}\ne U$, then  $U_{1}, \ldots , U_{t}\in
 \tau (G)$ by the choice of $U$ and so $U=U_{1}\times  \cdots \times  U_{t}\in
 \tau (G)$ since $\tau$ is special by hypothesis, a contradiction. Hence 
$U=U_{1}$ is a $\sigma _{i}$-group for some $i$.

(3) {\sl $U\leq \tau (M)$. In particular, $U \in \tau (M)$. }

Let $\pi =\pi (D)$. Then $G$ is $\pi$-separable by Condition (i)
 and $U$ is a $\pi'$-group since $D$ is a 
normal Hall $\pi$-subgroup of $G$ by hypothesis and $U\cap D=1$ by Claim (1).
 Then  $U^{d}\leq M$ for some $d\in D$,
 where $d^{-1}\in C_{G}(U^{d})$ by Lemma 2.6 and so  $U\leq M$.

 Therefore $U \in \tau (M)$ by hypothesis
 and so  $U \in \tau (G)$ since $\tau$ is special 
by hypothesis and $G$ is a $\sigma$-soluble $P\sigma T$-group by
 Conditions (i) and (ii)
 and Theorem 2.3. This contradiction completes the proof of the sifficiency
 of the condition   of the theorem.

Now assume that 
 every $\sigma$-subnormal subgroup of a $\sigma$-soluble group $G$
 belongs to the set $\tau (G). $  We show that in this case 
Conditions (i) and (ii) hold for  $G$. First note that 
 every $\sigma$-subnormal subgroup of $G$ $\sigma$-permutable in $G$ since 
$\tau$ is special by hypothesis. Therefore $G$ is a  $\sigma$-soluble $P\sigma T$-group
 and                                                    
so, in view of Theorem 2.3, $G=D\rtimes M$, where $ D=G^{\frak{N_{\sigma}}}$ 
and  Conditions (i) and (ii) in Theorem 2.3 hold for $D$ and $M$ hold. Therefore 
 we have only to show that  $H\in \tau (M) $ for every subgroup $H$ of $M$.

 Since
 $G/D$ is $\sigma$-nilpotent by Lemma 2.9, every subgroup $H/D\leq G/D$ is 
$\sigma$-subnormal  in $G/D$ by Lemma 2.10 and then, by Lemma 2.30(3) below,
 $H$ is $\sigma$-subnormal 
 in $G$ and so $H\in \tau (G)$. But then $H/D\in \tau (G/D)$.  Therefore for every 
subgroup $H/D$ of $G/D$  we have  $H/D\in \tau (G/D)$, so 
 for every 
subgroup $V$ of $M$  we have  $V\in \tau (M)$ since $M\simeq G/D$. Therefore Conditions 
(i) and (ii) hold for $G$.

The theorem is proved.

 We say that $G$ is a $Q\sigma T$-group  if every
 $\sigma $-subnormal subgroup of $G$ is   modular   in $G$.

It is clear that the lattice  of all subgroups ${\cal L}(G)$ of the group $G$  is modular 
if and only if every subgroup of $G$ is modular in $G$. Therefore, in view
 of Example 2.8(2), we get from Theorem 2.11 the following known result.

{\bf Corollary 2.12} (Hu,   Huang and Skiba  \cite{Hu12}). {\sl
A  group $G$ with
 $D=G^{\frak{N_{\sigma}}}$  is  a  $\sigma$-soluble $ Q\sigma T$-group if and only if 
 the following conditions hold:}

(i) {\sl $G=D\rtimes L$, where $D$   is an abelian  Hall
 subgroup of $G$ of odd order,  $L$ is  $\sigma$-nilpotent and 
  the lattice
of all subgroups ${\cal L}(L)$ of $L$ is modular,}
                                  
(ii) {\sl  every element of $G$ induces a
 power automorphism in $D$, and   }

(iii) {\sl  $O_{\sigma _{i}}(D)$ has
a normal complement in a Hall $\sigma _{i}$-subgroup of $G$ for all $i$.}

In view of Example 2.8(3),    we get from Theorem 2.11 in the case when 
 $\sigma =\sigma ^{1}=\{\{2\}, \{3\}, \{5\}, \ldots
\}$ the following classical result.

{\bf Corollart 2.13 } (Gasch\"{u}tz \cite{gasch}). {\sl  A group
  $G$ is a  soluble  $T$-group if and
 only if the following conditions are satisfied:}
 
(i) {\sl   the nilpotent residual $L$ of $G$ is an abelian Hall subgroup of odd order,}

(ii) {\sl  $G$ acts by conjugation on $L$ as a group power automorphisms, and }

(iii) {\sl  $G/L$ is
 a Dedekind  group}.

Every quasinormal subgroup is  clearly modular in the group. Moreover,
 the following remarkable fact is well-known.

{\bf Theorem 2.14} (Schmidt \cite[Theorem 5.1.1]{Schm}). {\sl A subgroup $A$ of $G$ is
 quasinormal in $G$ if and only if $A$ is modular and  subnormal in $G$}.

Recall that an \emph{Iwasawa group} is a group
 in which every subgroup is quasinormal.

In view of Example 2.8(2) and Theorem 2.14,    we get from Theorem 2.11
 in the case when 
 $\sigma =\sigma ^{1}=\{\{2\}, \{3\}, \{5\}, \ldots
\}$ the following well-known result.

{\bf Corollart 2.15 } (Zacher \cite{zaher}). {\sl  A group
  $G$ is a  soluble  $PT$-group if and
 only if the following conditions are satisfied:}
 
(i) {\sl   the nilpotent residual $L$ of $G$ is an abelian Hall
 subgroup of odd order,}

(ii) {\sl  $G$ acts by conjugation on $L$ as a group
 power automorphisms, and }

(iii) {\sl $G/L$ is an Iwasawa group}.

We say, following \cite{555}, that $G$ is a $T_{\sigma}$-group if every $\sigma$-subnormal
 subgroup of $G$  is normal.

In view of Example 2.8(3),    we get from Theorem 2.11 
 the following known result.

{\bf Corollary 2.16} (Zhang, Guo and Liu \cite{555}).
 {\sl   A    $\sigma$-soluble group $G$ with
 $D= G^{{\mathfrak{N}}_{\sigma}}$  is a  $T_{\sigma}$-group if  and only if 
   the following conditions hold:}

(i) {\sl $G=D\rtimes L$, where $D$   is an abelian  Hall
 subgroup of $G$ of odd order,  and $L$ is a  Dedekind group,}
                                  
(ii) {\sl  every element of $G$ induces a
 power automorphism in $D$, and   }

(iii) {\sl  $O_{\sigma _{i}}(D)$ has
a normal complement in a Hall $\sigma _{i}$-subgroup of $G$ for all $i$.}

{\bf Corollary 2.17} (Ballester-Bolinches, Pedraza-Aguilera and  P\`{e}rez-Calabuing
  \cite{Ad}). {\sl  A $\sigma$-soluble group $G$ is a $T_{\sigma}$-group if and only if 
$G$ is a $T$-group and the Hall $\sigma _{i}$-subgroups of $G$ are Dedekind
 for all $i\in I$.}

{\bf Proof. }  First assume that $G$ is a
 $\sigma$-soluble $T_{\sigma}$-group. Then $G$ satisfies Condition (i) and (ii) in  
Corollary 2.16, so $G$ is a $T$-group by Corollary 2.13. On the other hand, for every Hall
$\sigma _{i}$-subgroup $H$ of $G$ we have $H=(O_{i}(D))\times S$, where $O_{i}(D)$
 is a Hall subgroup of $D$ and $D$ is a Hall subgroup of $G$,  so $H$ is a  Dedekind group 
since $S\simeq DS/D$ and  $D$ are Dedekind.

Finally, suppose that  $G$ is a soluble $T$-group and the
 Hall $\sigma _{i}$-subgroups of $G$ are Dedekind
 for all $i\in I$. And let $D=G^{\mathfrak{N}_{\sigma}}$. Then $G/N$ is Dedekind by , so 
$D=G^{\mathfrak{N}}$. Hence   $G$ is a $\sigma$-soluble $T_{\sigma}$-group. 

The corollary is proved.

In the recent papers \cite{preprI, alg2023}, a description of  $P\sigma T$-groups $G$
 was obtained  for the case when every Hall $\sigma _{i}$-subgroup of
 $G$ is either supersoluble or $PST$-group for all $\sigma _{i}\in \sigma (G)$.

Our next goal to discuss some results of these two papers.

 {\bf Definition 2.18.}   We say    that
$(D, Z(D); U_{1},  \ldots , U_{k})$  is a 
\emph{Robinson $\sigma$-complex } (a \emph{Robinson complex} in the case where
 $\sigma =\sigma ^{1}=
\{\{2\}, \{3\}, \{5\},  \ldots \}$) of $G$ if $D\ne 1$ is a normal subgroup
 of $G$ such that: 

(i) $D/Z(D)=U_{1}/Z(D)\times \cdots \times U_{k}/Z(D)$, where $U_{i}/Z(D)$ is a  
 simple  non-$\sigma$-primary 
 chief factor of $G$, $Z(D)\leq \Phi(D)$, and

(ii) every chief factor of $G$  below $Z(D)$ is cyclic.

{\bf Example 2.19.}  (i) Let  $G=SL(2, 7)\times A_{7}\times A_{5}\times B$,
 where  $B=C_{43}\rtimes C_{7}$
is  a non-abelian group of order 301 and let 
 $\sigma =\{\{2, 3, 5\}, \{7, 43\}, \{2, 3, 5, 7, 43\}'\}$.  
  Then $$(SL(2, 7)\times A_{7}, Z(SL(2, 7)); SL(2, 7), A_{7}Z(SL(2, 7)))$$
 is a Robinson $\sigma $-complex of $G$ and
 $$(SL(2, 7)\times A_{7}\times A_{5}, Z(SL(2, 7)); SL(2, 7), A_{7}Z(SL(2, 7)),
 A_{5}Z(SL(2, 7)))$$  is a  Robinson complex of $G$.

(ii) If $(D, Z(D); U_{1},  \ldots , U_{k})$  is a 
Robinson $\sigma ^{\pi}$-complex of $G$ (see Example 2.1(ii)), then $U_{i}/Z(D)$ is
 neither 
a $\pi$-group nor a $\pi'$-group and we say in this case that
$(D, Z(D); U_{1},  \ldots , U_{k})$  is a 
Robinson \emph{$\pi, \pi'$-complex} of $G$.

(iii) If $(D, Z(D); U_{1},  \ldots , U_{k})$  is a 
Robinson $\sigma ^{1\pi}$-complex of $G$ (see Example 2.1(iii)), then $U_{i}/Z(D)$ is
neither
 a $\pi'$-group nor a $p$-group for all $p\in \pi$ and we say in this case that
$(D, Z(D); U_{1},  \ldots , U_{k})$  is a 
Robinson \emph{$\pi$-complex} of $G$.

 Let $\pi\subseteq \mathbb{P}$. If $\pi=\emptyset$, then we put $O_{\pi}(G)=O_{\emptyset}(G)=1$. We say that  \emph{$G$  satisfies
 ${\bf N}_{\pi}$}  if whenever $N$ is  a soluble normal
subgroup of $G$, $\pi'$-elements of $G$ induce power automorphisms in
  $O_{\pi}(G/N)$.  We also say,  following \cite[2.1.18]{prod}, that
  \emph{$G$  satisfies
 ${\bf N}_{p}$}  if whenever $N$ is  a soluble normal
subgroup of $G$, $p'$-elements of $G$ induce power automorphisms in
  $O_{p}(G/N)$.

Our next goal here is to prove the following fact.

{\bf Theorem 2.20.}  {\sl Suppose  that $G$ is a $\sigma$-full group  and every 
 Hall $\sigma _{i}$-subgroup 
 of $G$ is either  supersoluble or a $PST$-group for all $i\in I$.
 Then $G$ is a $P\sigma T$-group  if 
 and   only if  $G$  has a normal subgroup $D$ such that:}

(i) {\sl  $G/D$ is a $\sigma$-soluble $P\sigma T$-group,  }

(ii) {\sl if  $D\ne 1$,  $G$ has a Robinson $\sigma$-complex
 $(D, Z(D); U_{1},  \ldots , U_{k})$, and }

(iii) {\sl   for any set  $\{j_{1}, \ldots , j_{r}\}\subseteq \{1, \ldots , k\}$, where
 $1\leq r  < k$, $G$ and $G /U_{j_{1}}'\cdots U_{j_{r}}'$ satisfy
 ${\bf N}_{\sigma _{i}}$ for all $\sigma _{i}\in \sigma (Z(D))$.}

In view of Example  2.1(iii), we get from  Theorem 2.20 the following

{\bf Corollary 2.21.}  {\sl Suppose  that $G$ has a Hall $\pi'$-subgroup, where
$\pi= \{p_{1}, \ldots , p_{n}\}$, 
  and all such Hall subgroups of $G$
 are supersoluble. Then the condition
$\pi$-permutability is a transitive relation on $G$ if and only if
 $G$ has a normal  subgroup $D $ such that:}

(i) {\sl  $G/D$ is  $\pi $-soluble and the condition
$\pi$-permutability is a transitive relation on $G/D$, }

 (ii) {\sl if  $D\ne 1$, $G$ has a Robinson $\pi$-complex
 $(D, Z(D); U_{1},  \ldots , U_{k})$, and }

(iii)   {\sl for any set  $\{j_{1}, \ldots , j_{r}\}\subseteq \{1, \ldots , k\}$, where
 $1\leq r  < k$, $G$ and $G /U_{j_{1}}'\cdots U_{j_{r}}'$ satisfy
 ${\bf N}_{p}$ for all  $p\in \pi(Z(D))$ and, also,  
 ${\bf N}_{\pi'}$ for the case  $O_{\pi'}(Z(D))\ne 1$. }

 In view of Example  2.1(i), we get from Corollary   1.5 the following

{\bf Corollary   2.22} (Robinson  \cite{217}).  {\sl A group
  $G$ is an $PST$-group if
 and
only if  $G$ has a perfect normal  subgroup $D$ such that:}

(i) {\sl  $G/D$ is a soluble $PST$-group, }

(ii) {\sl if  $D\ne 1$,  $G$ has a Robinson complex
 $(D, Z(D); U_{1},  \ldots , U_{k})$, and }

(iii) {\sl   for any set  $\{j_{1}, \ldots , j_{r}\}\subseteq \{1, \ldots , k\}$, where
 $1\leq r  < k$, $G$ and $G /U_{j_{1}}'\cdots U_{j_{r}}'$ satisfy
 ${\bf N}_{p}$ for all $p\in \pi (Z(D))\cap \pi$.}

Theorem 2.20  has also many other consequences.
In particular, in view of Example  2.1(ii), we get from  Theorem 2.20 the following

{\bf Corollary 2.23.}  {\sl Suppose  that $G$ has a Hall $\pi$-subgroup
 and Hall $\pi'$-subgroup and all such Hall subgroups of $G$
 are supersoluble. Then the condition
$\pi, \pi'$-permutability is a transitive relation on $G$ if and only if
 $G$ has a normal  subgroup $D $ such that:}

(i) {\sl  $G/D$ is  $\pi $-separable and the condition
$\pi, \pi'$-permutability is a transitive relation on $G/D$, }

 (ii) {\sl if  $D\ne 1$, $G$ has a Robinson $\pi, \pi'$-complex
 $(D, Z(D); U_{1},  \ldots , U_{k})$, and }

(iii)   {\sl  for any set  $\{j_{1}, \ldots , j_{r}\}\subseteq \{1, \ldots , k\}$, where
 $1\leq r  < k$, $G$ and $G /U_{j_{1}}'\cdots U_{j_{r}}'$ satisfy
 ${\bf N}_{\pi}$ if $O_{\pi}(Z(D))\ne 1$
 and ${\bf N}_{\pi'}$ if $O_{\pi'}(Z(D))\ne 1$. }

{\bf Exapmle 2.24.} Let $\alpha: Z(SL(2, 5))\to Z(SL(2, 7))$ be an isomorphism and let 
$$D:= SL(2, 5) \Ydown SL(2, 7)=(SL(2, 5)\times SL(2, 7))/V,$$
 where $$V=\{(a, (a^{\alpha})^{-1})\mid a\in Z(SL(2, 5))\},$$
  is the direct product  of the groups $SL(2, 5)$ and $SL(2, 7)$ with joint center
 (see \cite[p. 49]{hupp}). 
Let   $M=(C_{23}\wr C_{11}) \Yup (C_{67}\rtimes C_{11}$) be
 the direct product  of the groups $C_{23}\wr C_{11}$ and $C_{67}\rtimes C_{11}$
 with joint 
factorgroup $C_{11}$  (see \cite[p. 50]{hupp}),
 where 
$C_{23}\wr C_{11}$ is the regular wreath product of the groups
$C_{23}$ and $ C_{11}$  and  $C_{67}\rtimes C_{11}$ is a non-abelian group of
 order 737.

 Now, let $G=D\times M$ and $\sigma =\{\{5\}, \{7\},  \{2, 11, 23\},  \{3, 67\},
  \{2, 3, 5, 7, 11, 23, 67\}'
\}.$  We show that $G$ is a 
$P\sigma T$-group. 
 In view of \cite[I, Satz 9.10]{hupp}, $D=U_{1}U_{2}$ and 
$U_{1}\cap U_{2}=Z(D)=\Phi (D)$, where  $U_{i}$ is normal in $D$, 
 $U_{1}/Z(D)$ is a a simple group of order 60, and 
$U_{2}/Z(D)$ is a a simple group of order 168. Hence $(D, Z(D); U_{1}, U_{2})$ is
 a Robinson    $\sigma$-complex   of $G$.  In view of \cite[I, Satz 9.11]{hupp},  $M$ has 
normal subgroups $R$ ($|R|=23^{11}$) and $L$ ($|L|=C_{67}$)
 such that $M/R\simeq C_{67}\rtimes C_{11}$ and $M/L\simeq  C_{23}\wr C_{11}.$ It
 is cleat that $M$ is not $\sigma$-nilpotent, so $M^{\mathfrak{N}_{\sigma}}=L$ since 
  $M/L$ is $\sigma$-primary, Therefore $M\simeq G/D$ is a
 $\sigma$-soluble $P\sigma T$-group by Theorem 2.3 and, clearly,
 $D=G^{\mathfrak{S}_{\sigma}}=G^{\mathfrak{S}}$.
  The group $G$ is $\sigma$-full and all Hall
 $\sigma _{i}$-subgroups of $G$ are supseroluble for all $i$. It is not also to 
show that $G$ satisfies ${\bf N}_{\pi}$, where $\pi=\{2, 11, 23\}$. Therefore
 Conditions
 (i), (ii), and (iii)  hold for $G$, so $G$ is a 
$P\sigma T$-group by Theorem 2.20.  

Assume thatt $G$ is a $PST$-group. Then, in view of Example 2.1(i) and Theorem 2.20,
 $M\simeq G/D$ is a soluble $PST$-group, so 
$M^{\mathfrak{N}}$  is a Hall subgroup of $M$ and all subgroups of 
$M^{\mathfrak{N}}$  are normal in $M$  by Example 2.1(i)   and Theorem 2.3.
 But $R$ contains subgroups
 which are 
not normal in $M$. Hence $M^{\mathfrak{N}}=L$, so  $M/L\simeq C_{23}\wr C_{11}$
 is nilpotent.  This contradiction shows that  $G$ is not a $PST$-group.

From Theorems 2.3 and 2.20 it follows that every $\sigma$-soluble $P\sigma T$-group is 
$\sigma$-supersoluble  and every $\sigma$-full $P\sigma T$-group with supersoluble Hall
$\sigma _{i}$-subgroups for all $i$ is a $\sigma$-$SC$-group
 in the sense of the following

{\bf Definition 2.25.}  We say   that $G$ is:

(1) \emph{$\sigma $-supersoluble}  \cite{3}
if every chief factor of $G$ below $G^{{\mathfrak{N}}_{\sigma}}$ is
cyclic;

(2) a  \emph{$\sigma $-$SC$-group}
if every chief factor of $G$ below $G^{{\mathfrak{N}}_{\sigma}}$ is  simple.

{\bf Example 2.26.}  (i)   $G$ is supersoluble if and only if it is $\sigma
$-supersoluble  where $\sigma =\sigma ^{1}$  (see Example 1.1(i)).

(ii) A group $G$ is called an \emph{$SC$-group} \cite{217} if every chief
factor of $G$ is a simple group.  Note that   $G$ is  an
$SC$-group   if and only if it is a $\sigma$-$SC$-group
 where $\sigma =\sigma ^{1}$.

(iii)     Let $G=A_{5}\times B$, where $A_{5}$ is
 the alternating group of degree 5 and $B=C_{29}\rtimes C_{7}$ is
 a non-abelian group of order 203, and let $\sigma =\{\{2, 3, 5\}, \{7\},
 \{29\}, \{2, 3, 5,  7, 29\}'\}$. Then  $G^{{\frak{N}}_{\sigma}}=C_{29}$,
 so
   $G$ is a $\sigma$-supersoluble group but it  is neither soluble
 nor  $\sigma$-nilpotent.

(iv)      Let $G=SL(2, 7)\times A_{7}\times A_{5}\times B$, where
  $B=C_{43}\rtimes C_{7}$
is  a non-abelian group of order 301 and let
 $\sigma =\{\{2, 3, 5\}, \{7, 43\}, \{2, 3, 5, 7, 43\}'\}$.
 Then  $G^{{\frak{N}}_{\sigma}}=SL(2, 7)\times A_{7}$, so
   $G$ is a $\sigma $-$SC$-group  but it  is not a $\sigma$-supersoluble group.

Let  $1\in \mathfrak{F}$ be a class of groups. Then $G^{\mathfrak{F}}$  is the
\emph{$ \mathfrak{F}$-residual} of $G$, that is, the intersection of all normal subgroups
 $N$ of $G$ with $G/N\in \mathfrak{F}$.
 The class  of groups $ 1\in\mathfrak{F}$ is
 said to be a \emph{formation} if every
 homomorphic image of $G/G^{\mathfrak{F}}$  belongs to $ \mathfrak{F}$ for every
 group $G$.
  The formation
$\mathfrak{F}$ is said to be \emph{(normally) hereditary }
 if $H\in \mathfrak{F}$ whenever $ G \in \mathfrak{F}$   and $H$ is a (normal) subgroup of $G$.

{\bf Lemma~2.27} (See  \cite[Proposition 2.2.8]{15}). {\sl Let $\frak{F}$ be a non-empty
formation and
$N$, $R$   subgroups of $G$, where $N$ is normal in $G$.}

(1) {\sl
 $(G/N)^{\frak{F}}=G^{\frak{F}}N/N.$  }

(2) {\sl If $G=RN$, then $G^{\frak{F}}N=R^{\frak{F}}N$}.

In what follows, 
    ${\mathfrak{U}}_{\sigma}$ is the class
 of all $\sigma $-supersoluble groups; ${\mathfrak{U}}_{c\sigma}$
  is the class of all $\sigma $-$SC$-groups.

In our proofs, we often use the following

{\bf Proposition~2.28.}    {\sl  For any partition $\sigma$ of  \ $\mathbb{P}$ the
following hold. }

(i) {\sl The class ${\mathfrak{U}}_{c\sigma}$ is a normally hereditary
 formation. }

(ii) {\sl The class ${\mathfrak{U}}_{\sigma}$ is a  hereditary
 formation } \cite{3}.

 {\bf Proof.}  (1) Let $D=G^{{\mathfrak{N}}_{\sigma}}$.
 First note that if
$R$ is a normal subgroup of $G$, then
$(G/R)^{{\mathfrak{N}}_{\sigma}}=DR/R$ by
Lemmas 2.9 and   2.27 and  so from the $G$-isomorphism
$DR/R\simeq D/(D\cap  R)$ we get that every chief factor of $G/R$ below
 $(G/R)^{{\mathfrak{N}}_{\sigma}}$ is simple if and only if every chief factor of $G$
between $D\cap R$ and $D$ is simple.
Therefore if  $G\in   {\mathfrak{U}}_{c\sigma}$, then $G/R\in   {\mathfrak{U}}_{c\sigma}$.
 Hence the class
${\mathfrak{U}}_{c\sigma}$  is closed under taking homomorphic images.

Now  we show that if  $G/R$, $G/N\in  {\mathfrak{U}}_{c\sigma}$,
 then $G/(R\cap N)  \in  {\mathfrak{U}}_{c\sigma}$.   We can assume without loss
of generality that $R\cap N=1$. Since $G/R \in  {\mathfrak{U}}_{c\sigma}$,
every chief factor of $G$
between $D\cap R$ and $D$ is simple. Also, every chief factor of $G$
between $D\cap N$ and $D$ is simple. Now let $H/K$ be any chief factor of $G$ below $D\cap R$.
Then $H\cap D\cap N=1$ and hence  from the $G$-isomorphism
 $$H(D\cap N)/K(D\cap N)\simeq H/(H\cap K(D\cap N))=H/K(H\cap D\cap N)=H/K$$
we get  that $H/K$ is simple since $D\cap N\leq K(D\cap N)\leq H(D\cap N) \leq D$. On the
other hand,    every chief factor of $G$
between $D\cap R$ and $D$ is also  simple.
 Therefore the Jordan-H\"{o}lder
theorem for groups with operators \cite[Ch. A, Theorem 3.2]{DH} implies that every
 chief factor of $G$
below $D$  is simple.  Hence   $G  \in
{\mathfrak{U}}_{c\sigma}$, so  the class
${\mathfrak{U}}_{c\sigma}$  is closed under taking subdirect products.

Finally, if  $H\trianglelefteq  G\in   {\mathfrak{U}}_{c\sigma}$, then from Lemmas 2.9
 and
2.27 and the isomorphism $H/(H\cap D)\simeq HD/D \in  {\mathfrak{N}}_{\sigma}$ we get that
$H^{{\mathfrak{N}}_{\sigma}}\leq  H\cap D$ and
so  every chief factor of $H$
below $H^{{\mathfrak{N}}_{\sigma}}$  is simple since every chief factor of $G$
below $D$ is simple. Hence
 $H\in
{\mathfrak{U}}_{c\sigma}$,  so  the class
${\mathfrak{U}}_{c\sigma}$  is closed under taking normal subgroups.

The proposition~is proved.

{\bf Proposition  2.29.}   {\sl Suppose  that  $G$ is a $P\sigma T$-group.
 Then }

(i) {\sl $G/R$  satisfies
 ${\bf N}_{\sigma _{i}}$ for every normal subgroup $R$ of $G$ and all $i\in I$, and }

(ii) {\sl if all
 Hall $\sigma _{i}$-subgroups of $G$ are supersoluble for all $i\in I$,
 then $G$ is a   $\sigma$-$SC$-group}.
 
To prove the proposition, we need a few lemmas.

Recall  that a subgroup $A$ of $G$ is called \emph{${\sigma}$-subnormal}
  in $G$ \cite{1} if   there is a subgroup chain  $$A=A_{0} \leq A_{1} \leq \cdots \leq
A_{n}=G$$  such that  either $A_{i-1} \trianglelefteq A_{i}$ or
$A_{i}/(A_{i-1})_{A_{i}}$ is  ${\sigma}$-primary
  for all $i=1, \ldots , n$.

We say that: an integer \emph{$n$ is a $\Pi$-number} if
 $\sigma (n)\subseteq \Pi$; a  subgroup $H$ of $G$ is a
 \emph{$\Pi$-subgroup   of $G$} if $|H|$ is a $\Pi$-number; 
  a $\sigma$-Hall  subgroup $H$ of $G$ is  a Hall
  \emph{$\Pi$-subgroup} of $G$  if    $H$ is  a $\Pi$-subgroup of $G$
 and $|G:H|$ is a $\Pi'$-number.    
We use  $O^{{\Pi}}(G) $ to denote the subgroup of $G$ generated by all  
its    ${\Pi}'$-subgroups.

{\bf Lemma~2.30.} {\sl Let  $A$,  $K$ and $N$ be subgroups of a $\sigma$-full group $G$.
 Suppose that   $A$
is $\sigma$-subnormal in $G$ and $N$ is normal in $G$.  }

(1) {\sl $A\cap K$    is  $\sigma$-subnormal in   $K$}.

(2) {\sl $AN/N$ is
$\sigma$-subnormal in $G/N$. }

(3) {\sl If $N\leq K$ and $K/N$ is
$\sigma$-subnormal in $G/N$, then $K$ is
$\sigma$-subnormal in $G.$}

(4) {\sl If $H\ne 1 $ is a Hall $\sigma _{i}$-subgroup of $G$  and $A$ is not  a
 $\sigma _{i}'$-group, then $A\cap H\ne 1$ is
 a Hall $\sigma _{i}$-subgroup of $A$. }

(5) {\sl If $A$ is a $\sigma _{i}$-group, then $A\leq O_{\sigma _{i}}(G)$.
}

(6) {\sl  If $A$ is a Hall $\sigma _{i}$-subgroup of $G$, then $A$ is normal in $G$.}

(7) {\sl If  $|G:A|$ is a $\Pi$-number,  then  $O^{\Pi}(A)=
 O^{\Pi}(G)$.}

 (8)    {\sl  If $O^{\sigma _{i}}(G)=G$ for all $i\in I$,
 then $A$ is subnormal in $G$. }

 (9)    {\sl  $A^{{\frak{N}}_{\sigma}}$    is subnormal in $G$. }

{\bf Proof. }    Assume that this Lemma~  is false and let $G$ be a counterexample of
 minimal     order.   By hypothesis, there is a subgroup chain  $A=A_{0} \leq
A_{1} \leq \cdots \leq A_{r}=G$ such that
either $A_{i-1} \trianglelefteq A_{i}$
  or $A_{i}/(A_{i-1})_{A_{i}}$ is  $\sigma $-primary  for all $i=1, \ldots , r$.
  Let   $M=A_{r-1}$.
  We can assume without loss of generality that $M\ne G$.

(1)--(7) See Lemma~2.6 in \cite{1}.

 (8) $A$ is
subnormal in $M$ by the choice of $G$. On the other hand, since $G$ is  $\sigma$-perfect,
 $G/M_{G}$  is not $\sigma$-primary. Hence $M$ is normal in $G$ and so $A$
is   subnormal in $G$.

(9)  $A$ is $\sigma$-subnormal in $AM_{G}\leq M$ by Part (1), so the
choice of $G$ implies that $A^{{\frak{N}}_{\sigma}}$ is
subnormal in $AM_{G}$.   Hence $G/M_{G}$ is a $\sigma
_{i}$-group for some $i$, so $M_{G}A/M_{G}\simeq A/(A\cap M_{G})$ is a $\sigma
_{i}$-group.  Hence $A^{{\frak{N}}_{\sigma}}\leq  M_{G}$, so $A^{{\frak{N}}_{\sigma}}$ is subnormal
 in $M_{G}$ and hence  $A^{{\frak{N}}_{\sigma}}$ is subnormal
 in $G$.

The lemma~is proved.

{\bf Lemma~2.31.}   {\sl  The following statements hold:}

(1) {\sl $G$ is a
 $P\sigma T$-group if and only if every  $\sigma$-subnormal subgroup of
$G$ is $\sigma$-permutable in $G$. }

(2) {\sl
If  $G$ is a
 $P\sigma T$-group, then every   quotient $G/N$ of $G$ is also a
$P\sigma T$-group. }

{\bf Proof.}   (1)  First note that if $A$ is a maximal 
 $\sigma$-subnormal  subgroup of $G$, then either $A$ is normal in $G$ or $G/A_{G}$ 
is a $\sigma _{i}$-group for some $i$. We show that $A$ is $\sigma$-permutable in $G$.
If $A$ is normal in $G$, it is clear. Now assume that $G/A_{G}$ 
is a $\sigma _{i}$-group. 
Let $H$ be a Hall $\sigma _{j}$-subroup of $G$. If $j\ne i$, then $H\leq A_{G}$ and so 
$AH=A=HA$. Finally, if $i=j$, then $A_{G}H=G$, so $AH=G=HA$.

Now  assume that $G$ is a $P\sigma T$-group and let $A$
 be a $\sigma$-subnormal subgroup of $G$. Then there is a subgroup chain 
$A=A_{0} \leq A_{1} \leq \cdots \leq
A_{n}=G$  such that   $A_{i-1}$ is a maximal 
 $\sigma$-subnormal  subgroup of $A_{i}$ and so  $A_{i-1}$ is $\sigma$-permutable in
 $A_{i}$  for all $i=1, \ldots , n$. But then $A$ is $\sigma$-permutable in
$G$.   Therefore every  $\sigma$-subnormal subgroup of any $P\sigma T$-group 
 is $\sigma$-permutable.

Finally, from Theorem B in \cite{1} it follow that   every  $\sigma$-permutable 
 subgroup
 of $G$ is  $\sigma$-subnormal in $G$.  Hence (1) holds.

(2)   Let $A/N$ be any $\sigma$-subnormal subgroup of $G/N$.
 Then $A$ is $\sigma$-subnormal $G$  by Lemma 2.30(3), so $A$
 is $\sigma$-permutable in $G$ and so $A/N$ is $\sigma$-permutable in $G$
 by Lemma 2.7(1).  Therefore we have (2) by Part (1).

The lemma is proved.

{\bf Lemma~2.32.}   {\sl
Let  $A$ and $B$ be subgroups of $G$, where $A$ is
 $\sigma$-permutable in  $G$. }

 (1) {\sl If $A\leq B$ and $B$ is $\sigma$-subnormal in $G$,
 then  $A$ is    $\sigma$-permutable  in $B$}.

(2) {\sl Suppose that $B$ is a $ \sigma _{i}$-group. Then $B$ is $\sigma$-permutable in 
 $G$ if and only if $O^{\sigma _{i}}(G) \leq N_{G}(B)
$} (See Lemma 3.1 in \cite{1}).

 {\bf Proof. }  (1) Let $\sigma (B) =\{\sigma _{1}, \ldots, \sigma _{n} \}$ and let 
$H_{i}$ be a Hall $ \sigma _{i}$-subgroup of $G$  for all $i$.  Let $x\in B$ and 
$H=H_{i}$. Then we have $AH^{x}=H^{x}A$,
 so $$AH^{x}\cap B=A(H^{x}\cap B)=A(H\cap B)^{x}=(H\cap B)^{x}A,$$ where 
$H\cap B$ is a  Hall $ \sigma _{i}$-subgroup of $B$ by Lemma~2.30(4).
  Hence  $A$ is
 $\sigma$-permutable in $B$ by Lemma 2.7(1).

The lemma~is proved.

 {\bf Proof of Proposition 2.29.}
  Let  $ S= G^{{\mathfrak{S}}_{\sigma}}$ be  the
 $\sigma$-soluble  residual and  $ D= G^{{\mathfrak{N}}_{\sigma}}$
 the $\sigma$-nilpotent residual of $G$.

(i)  In view of Lemma~2.31(2), we can assume without loss of
 generality that $R=1$.
 Let $L$ be   any soluble normal subgroup of $G$ and let $x$ be
 a $\sigma _{i}'$-element of $G$.  Let $V/L\leq O_{\sigma _{i}}(G/L)$.
Then $V/L$ is  $\sigma$-subnormal in $G/L$, so $V/L$ is
$\sigma$-permutable in $G/L$ by Lemma~2.7(1) since $G/L$ is a $P\sigma
T$-group by  Lemma~2.31(2). Therefore $xL\in O^{\sigma _{i}}(G/L)\leq N_{G/L}(V/L)$
 by Lemma~2.32(2).   Hence $G$  satisfies   ${\bf N}_{\sigma _{i}}$.

(ii) 
 Suppose that this  is false and let $G$ be a
counterexample of minimal order.  If $S=1$, then $G$ is $\sigma$-soluble
and so $G$ is a $\sigma$-$SC$-group  by Theorem B. Therefore $S\ne 1$, so $D\ne 1$.
Let $R$ be  a minimal normal subgroup of $G$ contained in $D$.
 Then $G/R$
is a $P\sigma T$-group by Lemma~2.31(2). Therefore the choice of $G$
implies that $G/R$ is a  $\sigma$-$SC$-group.    Since
 $(G/R)^{{\mathfrak{N}}_{\sigma }}=D/R$ by Lemmas 2.9 and 2.27,
every
chief  factor of $G/R$   below $D/R$ is simple. Hence
 every
chief  factor of $G$ between  $R$ and  $G^{{\mathfrak{N}}_{\sigma }}$ is
simple.
Therefore, in view of the Jordan-H\"{o}lder
theorem for groups with operators \cite[Ch. A, Theorem 3.2]{DH},
it is enough to show that   $R$ is   simple.

 Suppose
that this is  false and let $L$ be
a minimal normal subgroup of $R$.  Then $1  < L  < R$ and $L$
 is $\sigma$-subnormal in $G$,  so $L$ is $\sigma$-permutable in $G$
 by Lemma~2.31(1)   since  $G$ is  a $P\sigma T$-group. 
  Moreover, $L_{G}=1$ and so $L$  is
 $\sigma$-nilpotent by Theorem A.   Therefore  $R$ is a
$\sigma _{i}$-groups for some $i$ and so  $R\leq H$, where  $H$ is a 
Hall $\sigma _{i}$-subgroup of $G$. Since $H$ is supersoluble by hypothesis, $R$
 is abelian and so 
 there is a maximal subgroup $V$ of $R$ such that $V$ is normal in $H$.
 Therefore, in view of  Lemma~2.32(2),
we have   $G=HO^{\sigma _{i}}(G)\leq N_{G}(V).$  Hence $V=1$ and so $|R|=p$, a
contradiction. Thus $G$ is a $\sigma$-$SC$-group.

 The proposition is proved.

 {\bf Lemma~2.33.}   {\sl  Let $H/K$ be a non-abelian chief factor of $G$.
 If $H/K$ is simple, then $G/HC_{G}(H/K)$ is soluble.}

 {\bf Proof.}   Since   $C_{G}(H/K)/K=C_{G/K}(H/K)$, we can  assume without
loss of generality that $K=1$. Then  $G/C_{G}(H)\simeq V\leq
\text{Aut}(H)$ and $H/(H\cap C_{G}(H) )\simeq HC_{G}(H)/C_{G}(H)\simeq \text{Inn}(H)$
 since $C_{G}(H)\cap H=1.$ Hence
$$G/HC_{G}(H)\simeq  (G/C_{G}(H))/(HC_{G}(H)/C_{G}(H))\simeq W\leq
 \text{Aut}(H)/\text{Inn}(H).$$ From the validity of the Schreier
conjecture, it follows that $G/HC_{G}(H/K)$ is soluble.

The lemma~is proved.

{\bf Lemma~2.34.}   {\sl If $L$ is a non-abelian minimal subnormal subgroup of $G$, then
 $L^{G}$ is a minimal normal subgroup of $G$.}

{\bf Proof.}  Since every
 two perfect subnormal subgroups are permutable by \cite[II, Theorem 7.9]{26}, 
for some $x_{1}, \ldots , x_{t}\in G$ we have $L^{G}=L^{x_{1}}  \cdots  L^{x_{t}}$. 
Now let $R$ be a minimal normal subgroup of $G$ contained in $L^{G}$. In view of
 \cite[Ch. A, Lemma 14.3]{DH}, $R\leq N_{G}(L^{x_{i}})$, so
 for
 some $j$ we have  $L^{x_{j}}\leq R$ since $L^{G}$ and $R$ are non-abelian groups. 
But then  $L^{G}\leq R\leq L^{G}$, so  $L^{G}\leq R$.

The lemma~is proved.

  {\bf Theorem 2.35.}   {\sl Suppose that every $\sigma$-primary chief factor of $G$ 
is abelian. Then $G$ is a  $\sigma$-$SC$-group if and only if 
$G/G^{\mathfrak{S}_{\sigma}}$ is $\sigma$-supersoluble and 
 if  $G^{\mathfrak{S}_{\sigma}}\ne 1$, $G$ has a Robinson $\sigma$-complex
 $(G^{\mathfrak{S}_{\sigma}}, Z(G^{\mathfrak{S}_{\sigma}}); U_{1},  \ldots , U_{k}).$}

 {\bf Proof. } Let $D=G^{\mathfrak{S}_{\sigma}}$.  Then 
$O^{\sigma _{i}}(D)=D$ for all $i\in I$.

First assume that  $G$ is a  $\sigma$-$SC$-group. Then
 every chief factor of $G$  below $Z(D)$ is cyclic.    
  Now let $H/K$ be any chief factor of $G$ between $Z(D)$
 and $D$. If  $H/K$ is abelian, this factor  is cyclic, which implies that
 $D\leq C_{G}(H/K)$. On the other hand, if $H/K$ is a non-$\sigma$-primary
simple group, then Lemma~2.33 implies that $G/HC_{G}(H/K)$ is soluble.
Hence $$DHC_{G}(H/K)/HC_{G}(H/K)\simeq D/(D\cap HC_{G}(H/K))=D/HC_{D}(H/K)$$
is soluble, so  $D=HC_{D}(H/K)$ since $O^{\sigma _{i}}(D)=D$ for all $i$.  
 Therefore, in both cases,
every element of
$D$ induces an inner automorphism on $H/K$. Therefore $D$ is
quasinilpotent. Hence   in view of \cite[X, Theorem 13.6]{31},
 $D/Z(D)\simeq U_{1}/Z(D)\times \cdots \times U_{k}/Z(D)$, where $U_{j}/Z(D)$ is a
 non-$\sigma$-primary simple factor of $D$ for all $j=1, \ldots, k$.  
 Finally, note that $Z(D)\leq \Phi(D)$ since $O^{\sigma _{i}}(D)=D$ for all $i$. 
 Therefore  $(D,  Z(D); U_{1}, \ldots ,
U_{k})$ is a
Robinson $\sigma$-complex of  $G$ by Lemma~2.34 since $G$ is a $\sigma$-$SC$-group.

Now assume that $G/D$ is $\sigma$-supersoluble and, in the case $D\ne 1$,  
 $G$ has a Robinson $\sigma$-complex $(D, Z(D); U_{1},  \ldots , U_{k}).$ 
Then there is a chief series $1 =G_{0} < G_{1} < \cdots < G_{t-1} < G_{t}=D$ of $G$ 
 below $G^{\mathfrak{N}_{\sigma}} $  such that $ G_{i}/G_{i-1}$ is simple for all
 $i=1, \ldots, t$. 
Hence  the Jordan-H\"{o}lder
theorem for groups with operators \cite[Ch. A, Theorem 3.2]{DH}
 implies that every chief factor of $G$ below $G^{\mathfrak{N}_{\sigma}}$ is simple,
 that is, 
$G$ is a $\sigma$-$SC$-group.
 
 The theorem is proved.
                                                       
In the case $\sigma =\sigma ^{1}$ we get from Theorem 2.35 the following

{\bf Corollary 2.36}  (See Theorem 1.6.5 in \cite{prod}).  {\sl A group $G$ is a 
$SC$-group if and only if $G$ satisfies: }

(i) {\sl $G/G^{{\mathfrak{S}} }$ is supersoluble.}

 (ii) {\sl If $G^{\mathfrak{S}}\ne 1$, then  $G$ has a Robinson complex
 $(G^{\mathfrak{S}}, Z(G^{\mathfrak{S}}); U_{1},  \ldots , U_{k}).$}

{\bf Proposition 2.37.}  {\sl 
 If $G$ is a $\sigma$-$SC$-group with  $G^{{\mathfrak{S}} _{\sigma}}=
G^{{\mathfrak{U}} _{\sigma}}$, then
 $G^{{\mathfrak{S}} _{\sigma}}\leq
 C_{G}(G_{\mathfrak{S}}\cap G^{{\mathfrak{N}} _{\sigma}}).$  }
  
{\bf Proof.} Let $H/K$ be any chief factor of $G$ below 
$G_{\mathfrak{S}}\cap G^{{\mathfrak{N}} _{\sigma}}$. Then $H/K$ is cyclic since  
$G$ is a $\sigma$-$SC$-group, so $G_{\mathfrak{S}}\cap G^{{\mathfrak{N}} _{\sigma}}$ is
 contained in the supersoluble hypercentre $Z_{\mathfrak{U}}(G)$ of $G$. Then 
$G/C_{G}(G_{\mathfrak{S}}\cap G^{{\mathfrak{N}} _{\sigma}})$ is supersoluble by
 \cite[IV, Theorem 6.10]{DH} and so $G^{{\mathfrak{S}} _{\sigma}}=G^{{\mathfrak{U}} _{\sigma}}\leq 
  C_{G}(G_{\mathfrak{S}}\cap G^{{\mathfrak{N}} _{\sigma}})$.  
                     
 The proposition is proved.

In the case $\sigma =\sigma ^{1}$ we get from Proposition 2.37 the following

{\bf Corollary 2.38}  (See Proposition 1.6.4 in \cite{prod}).  {\sl If $G$ is
 an $SC$-group, then  $G^{\mathfrak{S}} \leq C_{G}(G_{\mathfrak{S}})$}.

{\bf Lemma~2.39.}   {\sl  Suppose that $G$ has a Robinson
 $\sigma$-complex $(D, Z(D); U_{1}, \ldots ,
U_{k})$ and let $N$ be a normal subgroup of $G$.  }

(1) {\sl If $N=U_{i}'$ and $k\ne 1$, then $Z(D/N)  =U_{i}/N =Z(D)N/N$ and 
 $$(D/N,
 Z(D/N); U_{1}N/N, \ldots , U_{i-1}N/N,
U_{i+1}N/N, \ldots , 
U_{k}N/N)$$ is  a Robinson $\sigma$-complex of $G/N$.  }

(2) {\sl If $N$ is nilpotent, then  $Z(DN/N)= Z(D)N/N$ and
 $$(DN/N, Z(DN/N); U_{1}N/N, \ldots ,
U_{k}N/N)$$ is  a Robinson $\sigma$-complex of $G/N$.}

 {\bf Proof. } Let $Z=Z(D)$. Then $U_{i}'Z=U_{i}$
 since $U_{i}/Z$ is a non-$\sigma$-primary simple group, so
 $U_{i}/U_{i}'=U_{i}'Z/U_{i}'\leq Z(D/U_{i}')$ and $U_{i}/U_{i}'\leq \Phi (D/U_{i}')$.

Moreover, 
  $DN/N$ is a non-identity normal subgroup of $G/N.$
Indeed,  if $D\leq N=U_{i}'\leq D$ for some $i$, then $D=U_{i}'=U_{i}=U_{1}$ and 
 so $k=1$,   
a contradiction. On the other hand,
 the case $D\leq N$, where $N$ is nilpotent, is also impossible since $U_{2}/Z$ is a
 non-$\sigma$-primary group.

 (1)  We can assume without loss of generality that $i=1.$

We have $$ D/U_{1}'=(U_{1}/U_{1}')(U_{2}U_{1}'/U_{1}')\cdots
 (U_{k}U_{1}'/U_{1}') ,$$
so $$(D/U_{1}')/(U_{1}/U_{1}')=((U_{2}U_{1}'/U_{1}')/(U_{1}/U_{1}'))\cdots
  ((U_{k}U_{1}'/U_{1}')/(U_{1}/U_{1}'))$$ and hence from $D\ne U_{1}$   and
 the $G$-isomorphisms  
 $$(U_{j}U_{1}'/U_{1}')/(U_{1}/U_{1}')\simeq U_{j}U_{1}'/U_{1}=U_{j}U_{1}/U_{1}\simeq 
U_{j}/(U_{1}\cap U_{j})=U_{j}/(ZU_{1}'\cap U_{j})$$$$=U_{j}/Z(U_{1}'\cap U_{j})=U_{j}/Z$$
we get that  $(D/U_{1}')/(U_{1}/U_{1}')$ is the  direct product of the non-$\sigma$-primary 
simple $G/U_{1}'$-invariant subgroups
 $(U_{j}U_{1}'/U_{1}')/(U_{1}/U_{1}')$, $j\ne 1$.    
Hence  $U_{1}/U_{1}'=ZU_{1}'/U_{1}'\leq Z(D/U_{1}') \leq U_{1}/U_{1}'$, so 
$U_{1}/U_{1}'=ZU_{1}'/U_{1}'=Z(D/U_{1}')\leq \Phi (D/U_{1}')$. 
Therefore   $$(D/U_{1}',
 Z(D/U_{1}'); U_{2}U_{1}'/U_{1}', \ldots , U_{k}U_{1}'/U_{1}')$$  
 is  a Robinson $\sigma$-complex of $G/N$.

(2) It is clear that $N\cap D\leq Z$, so  we have the $G$-isomorphisms 
$$ (DN/N)/(ZN/N)\simeq DN/ZN\simeq D/(D\cap NZ)=D/(D\cap N)Z=D/Z.$$  Note also that 
$$(U_{i}N/N)/(NZ/N)\simeq U_{i}/(NZ\cap U_{i})=U_{i}/(N\cap U_{i})Z=U_{i}/Z$$
 for all $i$ 
and 
$ (DN/N)/(ZN/N)$ is the direct product of the non-$\sigma$-primary simple 
 $G/N$-invariant  groups  $(U_{i}N/N)/(NZ/N)$.
  Hence $Z(DN/N)=ZN/N\leq \Phi (DN/N)$ 
 and 
every chief factor of $G/N$ below  $NZ/N\simeq _{G} Z/(N\cap Z) $ is cyclic. 
Therefore $$(DN/N, Z(DN/N); U_{1}N/N, \ldots ,
U_{k}N/N)$$ is  a Robinson $\sigma$-complex of $G/N$.

The lemma~is proved.

{\bf Lemma~2.40.}    {\sl Let $G$ be a non-$\sigma$-soluble $\sigma$-full
  $\sigma$-$SC$-group with  Robinson $\sigma$-complex
 $$(D, Z(D); U_{1}, \ldots ,
U_{k}),$$ where $D=G^{\mathfrak{S}_{\sigma}}=G^{\mathfrak{U}_{\sigma}}$.
  Let $U$ be a
 non-$\sigma$-permutable $\sigma$-subnormal subgroup of $G$ of minimal
order.   Then:}

(1) {\sl if $UU_{j}'/U_{j}'$ is $\sigma$-permutable  in  $G/U_{j}'$ for
all $j=1, \ldots, k$, then $U$ is $\sigma$-supersoluble, and }

(2) {\sl if $U$ is  $\sigma$-supersoluble and $UL/L$ is
$\sigma$-permutable  in  $G/L$ for
  all non-trivial nilpotent  normal subgroups $L$ of $G$, then
$U$ is a cyclic $p$-group for some prime $p$. }

{\bf Proof. }  Suppose that this Lemma~is false and let $G$ be a
counterexample of minimal order.  By hypothesis, for some
$i$ and for some  Hall $\sigma _{i}$-subgroup $H$ of $G$ we have $UH\ne
HU$. Moreover,  
 $O^{\sigma _{s}}(D)=D$ for all $s\in I$ and, in view of Proposition~2.28(ii),
 $G/D$ is $\sigma$-supersoluble.             

(1) Assume that this is false and suppose that
 $U\cap D\leq Z(D)$. Then every chief factor of $U$ below
 $U\cap Z(D)=U\cap D$ is cyclic  and, also,  $UD/D\simeq U/(U\cap  D)$ is
 $\sigma$-supersoluble by Proposition~2.28(ii).
 Hence  $U$ is $\sigma$-supersoluble, a contradiction. Therefore
 $U\cap D\nleq Z(D)$.   Moreover, Lemma~2.30(1)(2)
 implies that $(U\cap D)Z(D)/Z(D)$ is
$\sigma$-subnormal in $D/Z(D)$ and so  $(U\cap
D)Z(D)/Z(D)$ is a non-trivial
subnormal  subgroup of  $D/Z(D)$  by Lemma~2.30(8)
 since $O^{\sigma _{s}}(D/Z(D))=D/Z(D)$ for all $s$.

 Hence for some $j$ we
have $U_{j}/Z(D)\leq (U\cap
D)Z(D)/Z(D),$ so  $U_{j}\leq (U\cap
D)Z(D).$ But then $U_{j}'\leq  ((U\cap
D)Z(D))'\leq U\cap D.$  By hypothesis, $UU_{j}'/U_{j}'=U/U_{j}'$ is $\sigma$-permutable
  in
$G/U_{j}'$ and so
 $$UH/U_{j}'=(U
/U_{j}')(HU_{j}'/U_{j}')=(HU_{j}'/U_{j}')(U/U_{j}')=HU/U_{j}'.$$ Hence
 $UH=HU$,
 a contradiction. Therefore Statement (1) holds.

(2)  Let $N=
U^{{\mathfrak{N}}_{\sigma}}$.  Then $N$ is subnormal in $G$ by Lemma
2.30(9).
 Since $U$ is  $\sigma$-supersoluble by
hypothesis, $N < U$. By Lemmas 2.10,  2.27, and 2.30(3),  every proper  subgroup
 $V$ of $U$ with
$N\leq V$  is $\sigma$-subnormal in $G$, so the minimality of $U$ implies
that $VH=HV$. Therefore, if $U$ has at least two distinct  maximal subgroups $V$ and $W$
 such that $N\leq V\cap W$, then $U=\langle V, W \rangle $ is permutes with $H$ by
 \cite[Ch. A, Proposition 1.6]{DH},
 contrary to our assumption on $(U, H)$. Hence $U/N$
is a cyclic $p$-group for some prime $p$.
 Therefore we can  assume that $N\ne 1$.

 First assume that $p\in \sigma_{i} $. Lemma~2.30(4) implies that $H\cap U$
is a Hall $\sigma _{i}$-subgroup  of $U$, so $U=N(H\cap U)=(H\cap U)N$. Hence
$$UH=(H\cap U)NH=H(H\cap U)N=HU,$$ a contradiction. Thus
 $p\in \sigma_{j}$ for some $j\ne i$.

Now we show that $U$ is  a $P\sigma T$-group. Let $V$ be a proper
$\sigma$-subnormal subgroup of $U$.  Then $V$ is
$\sigma$-subnormal in $G$  since $U$ is
$\sigma$-subnormal in $G$.   The minimality of $U$ implies that $V$ is
  $\sigma$-permutable   in $G$,  so $V$ is
$\sigma$-permutable in $U$ by Lemma~ 2.32(1). Hence
 $U$ is a $\sigma$-soluble $P\sigma T$-group by Lemma~2.31(1), so $N$ is
an abelian Hall subgroup of  $U$ and all subgroups of $N$ are normal in $U$
 by Theorem B.
 Therefore $N$ is a $\sigma_{j}'$-group since $N=
U^{{\mathfrak{N}}_{\sigma}}$, so
 $N\leq O=O_{\sigma_{s}'}(F(G))$ by Lemma~2.30(5) (taking in the case
 $\sigma =\sigma ^{1}$).  By hypothesis, $OU/O$ permutes with
$OH/O$. By Lemma~2.30(1)(2), $OU/O$ is $\sigma$-subnormal in
$$(OU/O)(OH/O)=(OH/O)(OU/O)=OHU/O,$$
 where
$OU/O\simeq U/(U\cap O)$
 is a $\sigma _{j}$-group and $OH/O\simeq H/(H\cap O)$ is a $\sigma _{i}$-group.
Hence $UO/O$ is normal in  $OHU/O$ by Lemma~2.6(6). Therefore $H\leq N_{G}(OU).$

Now let $\Pi =\sigma \setminus \{\sigma _{j}\}$. 
Then 
 $$H\leq N_{G}(O^{\Pi}(OU))=N_{G}(O^{\Pi}(U))$$ by Lemma~2.30(7) since
 $|OU:U|=|O:O\cap U|$ is a $\Pi$-number and $U$ is $\sigma$-subnormal in $OU$ by Lemma 
2.30(1).
 For a Sylow $p$-subgroup $P$ of $U$ we have
$P\leq O^{\Pi}(U)$ since $p\in \sigma _{j}$. 
 Therefore $U=O^{\Pi}(U)N$ and $U/O^{\Pi}(U)\simeq 
N/(N\cap O^{\Pi}(U))$ is nilpotent and so $N\leq O^{\Pi}(U)$. But then $O^{\Pi}(U) = U$,
 so $H\leq N_{G}(U)$ and hence $HU=UH$, a
 contradiction. Therefore Statement  (2) holds.

The lemma~is proved.

{\bf Lemma~2.41} (See  Lemma~5 in 
\cite{knyag}). {\sl
Let $H$, $K$  and $N$ be pairwise permutable
subgroups of $G$ and  $H$  a Hall subgroup of $G$. Then $N\cap HK=(N\cap H)(N\cap K).$}

{\bf Lemma~2.42.}   {\sl  If $G$  satisfies
 ${\bf N}_{\sigma _{i}}$  and $N$ is a soluble normal  subgroup of $G$, then
$G/N$  satisfies
 ${\bf N}_{\sigma _{i}}$.}

 {\bf Proof. }  Let $L/N$ be a normal soluble subgroup of $G/N$,
 $(U/N)/(L/N)\leq      O_{\sigma _{i}}((G/N)/(L/N))$
    and   let $yN$ be a $\sigma _{i}'$-element in $G/N$. Then for some
 $\sigma _{i}'$-element $x\in G$    we have $yN=xN$.

 On the other hand, $L$ is a  soluble normal
  subgroup of $G$  and $U/L\leq      O_{\sigma _{i}}(G/L)$,
so $(U/L)^{x}=U^{x}/L=U/L$.  Hence
 $$((U/N)/(L/N))^{yN}=((U/N)/(L/N))^{xN}=(U/N)/(L/N).$$
Hence $G/N$  satisfies  ${\bf N}_{\sigma _{i}}$. 

The lemma~is proved.

{\bf Proof of Theorem 2.20.}  First assume that $G$ is a $P\sigma T$-groups and
 let $D=G^{{\mathfrak{S}}_{\sigma}}$  be the
 $\sigma$-soluble residual of $G$.  
 Then  $O^{\sigma _{i}}(D)=D$ for all $i\in I$.  
  Moreover, $G/D$ is a $\sigma$-soluble $P\sigma T$-group by Lemma~2.31(2), 
 hence   
$G/D$ is  $\sigma$-supersoluble by Theorem 2.3 and so, in fact, 
 $D=G^{{\mathfrak{U}}_{\sigma}}$  is the $\sigma$-supersoluble residual of $G$.  
 From Proposition 2.29(ii) it
 follows that $G$ is a $\sigma$-$SC$-group.

Therefore, if $D\ne 1$, then, in  view of Theorem 2.35,    
  $G$ has a Robinson $\sigma$-complex   
  $(D, Z(D);$ $U_{1},  \ldots , U_{k})$ 
 and, in view of Proposition 2.29(i),  
    for any set   $\{j_{1}, \ldots , j_{r}\}\subseteq \{1, \ldots , k\}$, where
 $1\leq r  < k$, $G$ and $G /U_{j_{1}}'\cdots U_{j_{r}}'$ satisfy
 ${\bf N}_{\sigma _{i}}$ for all $\sigma _{i}\in \sigma (Z(D))$.
  Therefore the necessity of the condition of the theorem holds.

 Now, assume, arguing 
by contradiction, that $G$ is non-$P\sigma T$-group of minimal order 
satisfying Conditions (i),  (ii) and (iii). 
 Then $D\ne 1$  
and, by Lemma~2.31(1),  $G$ has a $\sigma$-subnormal
subgroup $U$ such that $UH\ne HU$ for some $i\in I$ and some Hall $\sigma
_{i}$-subgroup $H$ of $G$ and, also, every $\sigma$-subnormal
subgroup $U_{0}$   of $G$   with $U_{0} < U$ is $\sigma$-permutable in
$G$.  From Conditions (i)  and (ii) it follows that 
$G$ is a $\sigma$-$SC$-group.

(1) {\sl   $U$ is $\sigma$-supersoluble.   }

 First assume that $k=1$, that is, $D=U_{1}=D'$.
 Then $Z(D)\leq \Phi (D)$
 and 
$(U\cap D)Z(D)/Z(D)$ is a $\sigma$-subnormal   subgroup 
 of a simple non-abelian group $D/Z(D)$ by Lemma~2.30(1)(2). Hence $U\cap D\leq Z(D)$,
 so $U\cap D=U\cap  Z(D)$.  Therefore
 every chief factor of $U$  below $U\cap D$ is cyclic. On the other hand, $U/(U\cap D)
\simeq UD/D$ is $\sigma$-supersoluble  by Condition (i), Theorem 2.3, and
 Proposition~2.28(ii)
 and so $U$ is $\sigma$-supersoluble.

Now let $k\ne 1$. We show that the hypothesis holds for $G/U_{j}'$
 for all $j=1, \ldots , k$. We can assume   without loss of generality that $j=1$.

  Let $N=U_{1}'$.  Then
  $ (G/N)/(D/N)\simeq G/D$ is a $\sigma$-soluble $P\sigma T$-group, so  Condition (i) 
holds for 
$G/N$.  From Lemma
 2.39(1) it follows that   
 $$(D/N, Z(D/N); U_{2}N/N, \ldots , U_{k}N/N)$$ is
 a Robinson $\sigma$-complex of $G/N$ and 
  $Z(D/N)=U_{1}/N=Z(D)N/N\simeq Z(D)/(Z(D)\cap N)$.
   Moreover, if
 $\{j_{1}, \ldots , j_{r}\}\subseteq \{2, \ldots , k\}$, where $2\leq r  < k$,
 then the quotients
 $G/N=G/U_{1}'$ and $$(G/N) /(U_{j_{1}}N/N)'\cdots (U_{j_{r}}N/N)'=
(G/N)/(U_{j_{1}}'\cdots U_{j_{r}}'U_{1}'/N)\simeq G/U_{j_{1}}'\cdots U_{j_{r}}'U_{1}'$$
 satisfy
 ${\bf N}_{\sigma _{i}}$ for all
 $\sigma _{i}\in  \sigma (Z(D/N))\subseteq \sigma (Z(D))$
 by Condition (iii). 

Therefore the hypothesis holds for $G/N=G/U_{1}'$, so the 
$\sigma$-subnormal subgroup $UU_{1}'/U_{1}'$ of $G/U_{1}'$  is $\sigma$-permutable in
  $G/U_{1}'$ by the choice of $G$. Finally,  
from Lemma 2.27(2) and Proposition 2.28(ii) 
 it follows that 
  $$(G/U_{1}')^{\mathfrak{S}_{\sigma}}=U_{1}'G^{\mathfrak{S}_{\sigma}}/U_{1}'
=DU_{1}'/U_{1}'=
 U_{1}'G^{\mathfrak{U}_{\sigma}}/U_{1}'= (G/U_{1}')^{\mathfrak{U}_{\sigma}}.$$
Therefore $U$ is $\sigma$-supersoluble by Lemma~2.40(1)

 (2) {\sl $U$ is a cyclic $p$-group for some prime $p\in \sigma _{s}$,
 where $s\ne i$.}

First we show that $U$ is a cyclic $p$-group for some prime $p$.

In view of Proposition 2.28(ii) and Lemma 2.31(2),
  $$(G/N)/(DN/N)\simeq G/DN\simeq (G/D)/(DN/D)$$ is
 a $\sigma$-soluble $P\sigma T$-group.  Moreover,  $Z(DN/N)= Z(D)N/N$ and  $(DN/N,  Z(D)N/N; U_{1}N/N, \ldots , $ $U_{k}N/N)$ is  a Robinson $\sigma$-complex of $G/N$ by Lemma~2.27(2).
In view of Lemma~2.42,  $G/N$   satisfies ${\bf N}_{\sigma _{i}}$ for all
$$\sigma _{i}\in \sigma (Z(D)N/N)\subseteq \sigma (Z(D)).$$ 
Similarly, for any set 
 $\{j_{1}, \ldots , j_{r}\}\subseteq \{1, \ldots , k\}$, where
 $1\leq r  < k$,  $$(G/N) /(U_{j_{1}}N/N)'\cdots (U_{j_{r}}N/N)'=
(G/N) /(U_{j_{1}}'\cdots U_{j_{r}}'N/N)$$$$\simeq G/U_{j_{1}}'\cdots U_{j_{r}}'N\simeq
   (G/U_{j_{1}}'\cdots U_{j_{r}}')/(U_{j_{1}}'\cdots
U_{j_{r}}'N/U_{j_{1}}'\cdots U_{j_{r}}')
$$  satisfies ${\bf N}_{\sigma _{i}}$ for all
$\sigma _{i}\in \sigma (Z(D)N/N)$ by Lemma~2.42 since $$U_{j_{1}}'\cdots
U_{j_{r}}'N/U_{j_{1}}'\cdots U_{j_{r}}'\simeq N/(N\cap U_{j_{1}}'\cdots
U_{j_{r}}')$$ is nilpotent.
 Therefore   the hypothesis holds
on $G/N$, so $UN/N$ is $\sigma$-permutable in  $G/N$. 
Also we have $(G/N)^{\mathfrak{S}_{\sigma}}= (G/N)^{\mathfrak{U}_{\sigma}}.$
 Therefore    $U$ is a cyclic $p$-group for some prime $p\in \sigma _{s}$ by Lemma
 2.40(2).  
 Finally, Lemma
 2.30(4)
implies that in the case $i=s$ we have $U\leq H$, so $UH=H=HU$. Therefore
 $s\ne i$.

(3) $\sigma _{s}\not \in \sigma (Z(D))$ (This follows from Condition (iii)
and  Lemma~2.32(2) since
 $H$ is a $\sigma _{s}'$-group by Claim (2)).

 (4) $O_{\sigma _{s}}(G)\cap D=1$.

Assume $O_{\sigma _{s}}(G)\cap D\ne 1$. Since  
 $D/Z(D)$ is  the direct product of non-$\sigma$-primary
 simple groups, $O_{\sigma _{s}}(G)\cap D\leq
Z(D)$. But then $\sigma _{s}
\in \sigma (Z(D))$, contrary to Claim (3).
 Therefore we have (4).

{\sl Final contradiction.}  By Lemma~2.30(2), $UD/D$ is
$\sigma$-subnormal in $G/D$. On the other hand, $HD/D$ is a Hall $\sigma
_{i}$-subgroup of $G/D$. Hence $$(UD/D)(HD/D)=(HD/D)(UD/D)=HUD/D$$ by Condition (ii) and
 Lemma~2.31(1),
 so $HUD$ is a subgroup of $G$.
                  Therefore,   by Claim  (4) and Lemma~2.41,
 $$UHD\cap HO_{\sigma _{s}}(G)
=UH(D\cap HO_{\sigma _{s}}(G))=UH(D\cap H)(D\cap
O_{\sigma _{s}}(G))=UH(D\cap H)=UH$$ is a subgroup of $G$ and so $HU=UH$,
 a contradiction.  Therefore the sufficiency of the condition of the theorem holds.

The theorem is proved.

\section{Groups in which every $\sigma$-subnormal  subgroup is $\sigma$-quasinormal}

The quasinormal and the Sylow permutable subgroups  have many useful for applications
  properties. 
For instance, if $A$ is quasinormal in $G$, then:  {\sl $A$ is subnormal in $G$}
 (Ore \cite{5}),  {\sl $A/A_{G}$ is nilpotent} (Ito and Szep \cite{It}), 
 {\sl every chief factor 
$H/K$ of $G$ between $A_{G}$ and $A^{G}$ is central, that is,  $C_{G}(H/K)=G$}
 ({Maier and Schmid  \cite{MaierS}),
  and,
 in general,  the section \emph{$A/A_{G}$ is not necessarily abelian}
(Thomson  \cite{Th}).

Note also that the quasinormal  subgroups have a 
 close connection with the  modular subgroups.

Every quasinormal is clearly modular in the group. Moreover,
 the following remarkable fact is well-known.

{\bf Theorem 3.1} (Schmidt \cite[Theorem 5.1.1]{Schm}). {\sl A subgroup $A$ of $G$ is
 quasinormal in $G$ if and only if $A$ is modular and  subnormal in $G$}.

 This  result made
 it possible to find analogues of quasinormal subgroups in the theory of the
 $\sigma$-properties of a group.

{\bf Definition 3.2.} We say, following \cite{Hu11}, that a subgroup $A$ of $G$
 is \emph{$\sigma$-quasinormal} in $G$ if  $A$ is modular and  $\sigma$-subnormal in $G$.

The   description of the 
 $PT$-groups
 was first obtained by   Zacher \cite{zaher},
  for the soluble  case, and
  by Robinson in \cite{217}, for the general case.
In the   further publications,   authors
 (see Chapter 2 in \cite{prod}) have found out and  described
many other   interesting characterizations   of  
 $PT$-groups and generalized $PT$-groups.

The
 theory of $\sigma$-quasinormal subgroups was constructed in the paper \cite{Hu11}.
 In particular,  it is proved the following result covering in the
 case $\sigma =
\sigma ^{1}=\{\{2\}, \{3\}, \{5\} \ldots  \}$ the above mentioned results in
 \cite{5, It, MaierS}.

{\bf Theorem 3.3} (See Theorem C in \cite{Hu11}). {\sl 
 Let $A$ be a $\sigma$-quasinormal subgroup of   $G$.
 Then the following statements hold:}

(i) {\sl If $G$
possesses a  Hall $\sigma_{i}$-subgroup,
 then $A$ permutes with each Hall $\sigma_{i}$-subgroup  of  $G$. }

(ii) {\sl The quotients
$A^{G}/A_{G}$ and $G/C_{G}(A^{G}/A_{G}) $ are $\sigma$-nilpotent, and }

(iii) {\sl Every chief factor of $G$ between
 $A^{G}$ and $A_{G}$ is $\sigma$-central in $G$.  }

(iv)  {\sl For every $i$ such that $\sigma _{i} \in \sigma
(G/C_{G}(A^{G}/A_{G}))$  we have
 $\sigma _{i} \in  \sigma (A^{G}/A_{G}).$
}

(v) {\sl $A$  is $\sigma$-seminormal in $G$.}

However, the following problem still remains open.

{\bf Question 3.4.}  {\sl What is the structure of  the $Q\sigma T$-groups,
 that are, groups in which every $\sigma$-subnormal  subgroup is $\sigma$-quasinormal?}

Partially,  Problem 3.4 was solved in the recent paper \cite{Hu12}.

{\bf Theorem 3.5 } (See Theorem 1.5 in \cite{Hu12}). {\sl
  Let $D$ be the $\sigma$-nilpotent residual of $G$, that is, the
 intersection of all normal subgroups $N$
 of $G$ with $\sigma$-nilpotent quotient $G/N$.
 If  $G$ is  a  $\sigma$-soluble  $Q\sigma T$-group,
 then   the following conditions hold:}

(i) {\sl $G=D\rtimes L$, where $D$   is an abelian  Hall
 subgroup of $G$ of odd order,  $L$ is  $\sigma$-nilpotent and 
  the lattice
of all subgroups ${\cal L}(L)$ of $L$ is modular,}
                                  
(ii) {\sl  every element of $G$ induces a
 power automorphism in $D$, and   }

(iii) {\sl  $O_{\sigma _{i}}(D)$ has
a normal complement in a Hall $\sigma _{i}$-subgroup of $G$ for all $i$.}

{\sl Conversely, if  Conditions (i), (ii)  and (iii) hold for  some
 subgroups $D$ and $L$ of
 $G$, then every $G$ is a $\sigma$-soluble $Q\sigma T$-group.}

In the recent paper  \cite{preprI}, Problem  3.4
 was solved in the general case.

 Let $\pi \subseteq \Bbb{P}$. Then we say \emph{$G$  satisfies}  
 ${\bf P}_{\pi}$  if whenever $N$ is  a soluble normal
subgroup of $G$ and $G$ has a Hall
 $\pi$-subgroup, every subgroup of $O_{\pi}(G/N)$ is modular
 in   Hall   $\sigma _{i}$-subgroups of $G/N$.

{\bf Theorem 3.6}  (Safonova, Skiba \cite{preprI})
 {\sl Suppose  that $G$ is a $\sigma$-full group.
 Then $G$ is a $Q\sigma T$-group  if 
 and   only if  $G$  has a normal subgroup $D$ such that:}

(i) {\sl  $G/D$ is a $\sigma$-soluble $Q\sigma T$-group, }

(ii) {\sl if  $D\ne 1$,  $G$ has a Robinson $\sigma$-complex
 $(D, Z(D); U_{1},  \ldots , U_{k})$ and }

(iii) {\sl   for any set $\{i_{1}, \ldots , i_{r}\}\subseteq \{1, \ldots , k\}$, where
 $1\leq r  < k$,  $G$ and $G /U_{i_{1}}'\cdots U_{i_{r}}'$ satisfy
 ${\bf N}_{\sigma _{i}}$ for all $\sigma _{i}\in \sigma (Z(D))$ and
 ${\bf P}_{\sigma _{i}}$ 
for all $\sigma _{i}\in \sigma (D)$. }

In the case $\sigma =\sigma ^{1\pi}$
 (see Example 2.1(ii)) we get from Theorem 3.6 the following

{\bf Corollary 3.7.}  {\sl Suppose that  $G$ has a Hall $\pi'$-subgroup.
 Then every  $\pi$-subnormal subgroup of $G$ is modular in $G$  if 
 and   only if  $G$  has a normal subgroup $D$ such that:}

(i) {\sl  $G/D$ is a $\pi$-soluble group in which every $\pi$-subnormal
 subgroup  is modular,  }

(ii) {\sl if  $D\ne 1$,  $G$ has a Robinson $\pi$-complex
 $(D, Z(D); U_{1},  \ldots , U_{k})$ and  }

(iii) {\sl   for any set  $\{i_{1}, \ldots , i_{r}\}\subseteq \{1, \ldots , k\}$, where
 $1\leq r  < k$,    
$G$ and $G /U_{i_{1}}'\cdots U_{i_{r}}'$ satisfy   ${\bf N}_{p}$
 for all  primes $p$ dividing $|Z(D)|$ and 
 ${\bf N}_{\pi'}$ if  $O_{\pi'}(Z(D))\ne 1$, and, also  
 $G$ and $G /U_{i_{1}}'\cdots U_{i_{r}}'$ satisfy  ${\bf P}_{p}$ for all
  primes $p$ dividing $|D|$ and    ${\bf P}_{\pi'}$ 
 if  $O_{\pi'}(D)\ne 1$.}

In the case  $\pi=\Bbb{P}$, we get from Corollary 3.7 the following  classical result.

{\bf Corollary 3.8} (Robinson \cite{217}). {\sl $G$ is a $PT$-group if 
 and   only if  $G$  has a normal perfect subgroup $D$ such that:}

(i) {\sl  $G/D$ is a soluble $PT$-group, and }

(i) {\sl if $D\ne 1$, $G$ has a Robinson complex
 $(D, Z(D); U_{1},  \ldots , U_{k})$ and }

(iii) {\sl   for any set  $\{i_{1}, \ldots , i_{r}\}\subseteq \{1, \ldots , k\}$, where
 $1\leq r  < k$,  $G$ and $G /U_{i_{1}}'\cdots U_{i_{r}}'$ satisfy
 ${\bf N}_{p}$ for all $p\in \pi (Z(D))$ and
 ${\bf P}_{p}$ 
for all $p \in \pi (D)$. }

Theorem 3.6 has also many other consequences.
In particular, in view of Example  2.1(ii), we get from  Theorem 3.6 the following

{\bf Corollary 3.9.}  {\sl Suppose that  $G$ has  a Hall
 $\pi$-subgroup and a Hall $\pi'$-subgroup.
 Then every  $\pi, \pi'$-subnormal subgroup of $G$ is modular in $G$  if 
 and   only if  $G$  has a normal subgroup $D$ such that:}

(i) {\sl  $G/D$ is a $\pi$-separable group in which every $\pi, \pi'$-subnormal
 subgroup is modular,  }

(ii) {\sl if  $D\ne 1$,  $G$ has a Robinson $\pi, \pi'$-complex
 $(D, Z(D); U_{1},  \ldots , U_{k})$ and  }

(iii) {\sl   for any set  $\{i_{1}, \ldots , i_{r}\}\subseteq \{1, \ldots , k\}$, where
 $1\leq r  < k$,    
$G$ and $G /U_{i_{1}}'\cdots U_{i_{r}}'$ satisfy   ${\bf N}_{\pi}$ if
  $O_{\pi}(Z(D))\ne 1$
and 
 ${\bf N}_{\pi'}$ if  $O_{\pi'}(D)\ne 1$, and, also,  
 $G$ and $G /U_{i_{1}}'\cdots U_{i_{r}}'$ satisfy    ${\bf P}_{\pi}$ 
 if  $O_{\pi}(D)\ne 1$ and  ${\bf P}_{\pi'}$ 
 if  $O_{\pi'}(D)\ne 1$.}

{\bf Exapmle 3.10.} Let $\alpha: Z(SL(2, 5))\to Z(SL(2, 7))$ be an isomorphism and let 
$$D:= SL(2, 5) \Ydown SL(2, 7)=(SL(2, 5)\times SL(2, 7))/,$$
 where $$V=\{(a, (a^{\alpha})^{-1})\mid a\in Z(SL(2, 5))\},$$
  is the direct product  of the groups $SL(2, 5)$ and $SL(2, 7)$ with joint center
 (see \cite[p. 49]{hupp}). 
Let   $M=(C_{23}\wr C_{11}) \Yup
 (C_{67}\rtimes C_{11}$) (see \cite[p. 50]{hupp}, where 
$C_{23}\wr C_{11}$ is the regular wreath product of the groups
$C_{23}$ and $ C_{11}$  and  $C_{11}\wr C_{67}$ is a non-abelian group of
 order 737. Then $M$ is a soluble $P\sigma T$-group and $M$ is not  a $PST$ by
 Theorem B. 
 Let  $G=D\times M$ and  $\sigma =\{\{2, 3\}, \{5, 11, 23\},
 \{7, 67\},
 \{2, 3, 5, 7,   11, 23, 67\}'\}$. Then $G$ is not a $PT$-group by Corollary 3.8.
 Moreover, $G$
 is $\sigma$-ful  and
Conditions (i), (ii), and (iii) in Theorem 3.6 hold for $G$. Therefore 
every  $\sigma$-subnormal subgroup of $G$ is modular in $G$.

To prove this theorem, we need the results from Section 2 and several new lemmas.

First note that from Theorem 3.3(1) we get the following  

{\bf Lemma  3.11.}  {\sl  Suppose that   $G$  is a $Q\sigma T$-group. Then }

(1) {\sl  every  $\sigma$-subnormal subgroup of $G$ is $\sigma$-seminormal in $G$, and }
 In particular,

(2) {\sl if $G$ is $\sigma$-full, then every
 $\sigma$-quasinormal subgroup of $G$ is $\sigma$-permutable in $G$.  
}

From Lemma  3.11 we get the following

{\bf Corollary 3.12.}  {\sl Every $\sigma$-full $Q\sigma T$-group is  a
 $P\sigma T$-group. }

{\bf Lemma 3.13.}   {\sl  If  $G$ is an
 $Q\sigma T$-group, then every   quotient $G/N$ of $G$ is also an
$Q\sigma T$-group. }

{\bf Proof.}  
Let $L/N$ is a $\sigma$-subnormal subgroup of $G/N$. Then $L$ is 
$\sigma$-subnormal subgroup in $G$ by Lemma 2.30(3), so $L$
 is modular in $G$   and
 then $L/N$ is modular and so $\sigma$-quasinormal
 in $G/N$ by  \cite[Page 201, Property (3)]{Schm}. Hence 
$G/N$  is  a  $Q\sigma T$-group.

 The lemma is proved.

{\bf Lemma 3.14.}   {\sl Suppose  that $G$ is a $\sigma$-full  $Q\sigma T$-group. Then 
$G/R$  satisfies   ${\bf P}_{\sigma _{i}}$  and ${\bf N}_{\sigma _{i}}$ 
  for every normal subgroup $R$ of $G$ and all $i\in I$.}

{\bf Proof.}  
 In view of Lemma 3.12, we can assume without loss of
 generality that $R=1$. Let  $U/N\leq O_{\sigma _{i}}(G/N)$. Then $U/N$ is
 $\sigma$-subnormal in $G/N$, so $U$ is   $\sigma$-subnormal in $G$ by Lemma 2.30(3).
 Thereore
$U$ is modular in $G$ 
since $G$ is a $Q\sigma T$-group and 
 so, by  \cite[Page 201, Property (3)]{Schm},  $U/R$ is
 modular in every  Hall $\sigma _{i}$-subgroup $H/R$ of $G/N$ since 
$O_{\sigma _{i}}(G/N)\leq H/N$. 
 Hence $G$  satisfies 
 ${\bf P}_{\sigma _{i}}$.  Moreover, $U/N$ is
 $\sigma$-permutable in $G/N$ by Lemma 3.11 and  so for every $\sigma _{i}'$-element 
$x\in G$ we have  $xL\in O^{\sigma _{i}}(G/N)\leq N_{G/L}(U/N)$
 by Lemma 2.32(2).   Hence $G$  satisfies   ${\bf N}_{\sigma _{i}}$.

The lemma is proved.

{\bf Lemma 3.15.}  {\sl  Let $G$ be a non-$\sigma$-soluble $\sigma$-full
  $\sigma$-$SC$-group with  Robinson $\sigma$-complex
 $$(D, Z(D); U_{1}, \ldots ,
U_{k}),$$ where $D=G^{\mathfrak{S}_{\sigma}}=G^{\mathfrak{U}_{\sigma}}$.
  Let $U$ be a  $\sigma$-subnormal  non-modular  (non-normal)
  subgroup of $G$ of minimal
order.   Then:}

(1) {\sl If $UU_{i}'/U_{i}'$ is modular (respectively, normal) in  $G/U_{i}'$ for
all $i=1, \ldots, k$, then $U$ is $\sigma$-supersoluble.}

(2) {\sl If $U$ is  $\sigma$-supersoluble and $UL/L$ is modular (respectively, normal)
 in  $G/L$ for
  all non-trivial nilpotent  normal subgroups $L$ of $G$, then
$U$ is a cyclic $p$-group for some prime $p$. }

{\bf Proof. }  Suppose that this lemma is false and let $G$ be a 
counterexample of minimal order.

(1)  Assume this is false. Suppose that 
$U\cap D\leq Z(D)$. Then every chief factor of $U$ below
 $U\cap Z(D)=U\cap D$ is cyclic  and, also,  $UD/D\simeq U/(U\cap  D)$ is
 $\sigma$-supersoluble by Proposition~2.28.
 Hence  $U$ is $\sigma$-supersoluble, a contradiction. Therefore
 $U\cap D\nleq Z(D)$.  
$U\cap D\nleq Z(D)$.   Moreover, Lemma~2.30(1)(2)
 implies that $(U\cap D)Z(D)/Z(D)$ is
$\sigma$-subnormal in $D/Z(D)$ and so  $(U\cap
D)Z(D)/Z(D)$ is a non-trivial
subnormal  subgroup of  $D/Z(D)$  by Lemma~2.30(8)
 since $O^{\sigma _{s}}(D/Z(D))=D/Z(D)$ for all $s$.

 Hence for some $j$ we
have $U_{j}/Z(D)\leq (U\cap
D)Z(D)/Z(D),$ so  $U_{j}\leq (U\cap
D)Z(D).$ But then $U_{j}'\leq  ((U\cap
D)Z(D))'\leq U\cap D.$  By hypothesis, $UU_{j}'/U_{j}'=U/U_{j}'$ is modular
 (respectively, normal) in 
  in
$G/U_{j}'$ and so  $U$ is modular (respectively, normal) in $G$ by \cite[p.~201, Property~(4)]{Schm},
 a contradiction. Therefore Statement (1) holds.

(2)  Let $N=
U^{{\mathfrak{N}}_{\sigma}}$.  Then $N$ is subnormal in $G$ by Lemma  2.30(9).
 Since $U$ is  $\sigma$-supersoluble by
hypothesis, $N < U$. By Lemmas 2.9,  2.10, and 2.30(3),
  every proper  subgroup
 $V$ of $U$ with
$N\leq V$  is $\sigma$-subnormal in $G$, so the minimality of $U$ implies
that $V$ is modular (respectively, normal) in $G$.

 Therefore, if $U$ has at least two distinct  maximal subgroups $V$ and $W$
 such that $N\leq V\cap W$, then $U=\langle V, W \rangle $ 
 is modular  (respectively, normal) in 
$G$  by \cite[p. 201, Property (5)]{Schm}, contrary to our assumption on $U$.
 Hence $U/N$ 
is a cyclic $p$-group for some $p\in \sigma _{i}$ and $N\ne 1$ since $U$ is not cyclic.

Now we show that $U$ is  a $P\sigma T$-group. Let $V$ be a proper 
$\sigma$-subnormal subgroup of $U$.  Then $V$ is  
$\sigma$-subnormal  in $G$  since $U$ is  
$\sigma$-subnormal in $G$, so $V$ is modular (respectively, normal)
 in $G$ and hence $V$ is 
modular (respectively, normal) in $U$.
 Therefore $V$ is $\sigma$-permutable in $U$ by Lemma 3.11. Hence 
 $U$ is a $\sigma$-soluble soluble $P\sigma T$-group, so $N=U^{{\mathfrak{N}_{\sigma}}}$
 is a 
Hall abelian   subgroup of $U$ and every subgroup of $N$ is normal in $U$ by
Theorem 2.3. Therefore, in fact, for a Sylow $p$-subgroup $P$ of $U$
 we have $U=N\rtimes P$ and $P$ is a cyclic Hall $\sigma _{i}$-subgroup of $U$. Let

 Clearly,  
$N$ is $\sigma$-quasinormal  in $G$.  Assume that for some minimal normal
 subgroup $R$ of $G$ we have  $R\leq N_{G}$. Then $R$ is abelian, 
 $U/R$ is modular (respectively, normal) in $G/R$  by hypothesis,
 so  $U$ is modular (respectively, normal) in $G$ by
 \cite[p. 201, Property (4)]{Schm},
 a contradiction. Therefore 
$N_{G}=1$, so $P\leq C_{G}(N^{G})$ since for every $k$ such that $\sigma _{k} \in \sigma
(G/C_{G}(N^{G})$  we have
 $\sigma _{i} \in  \sigma (N^{G})$  by Theorem 3.3(iv).
 Therefore $U=N\times P$ is
 $\sigma$-nilpotent and so $N=1$, a contradiction. 
 Therefore Statement  (2) holds.     

The lemma is proved.

{\bf Lemma 3.16.}   {\sl Let $N$ be a soluble normal  subgroup of $G$.}
   {\sl  If $G$  satisfies   $P_{\sigma _{i}}$, then $G/N$  satisfies
 $P_{\sigma _{i}}$.}

 {\bf Proof. }  Let $L/N$ be a normal soluble subgroup of $G/N$ and 
$(U/N)/(L/N)\leq      O_{\sigma _{i}}((G/N)/(L/N))$. Then
 $L$ is a  soluble normal
  subgroup of $G$  and $U/L\leq O_{\sigma _{i}}(G/L)$.

  By hypothesis, $U/L$ is modular in a Hall $\sigma _{i}$-subgroup $H/L$ of $G/L$. 
Then  $(U/N)/(L/N)$ is modular in the Hall $\sigma _{i}$-subgroup 
 $(H/N)/(L/N)$ of $(G/N)/(L/N)$. 
Hence $G/N$  satisfies  $P_{\sigma _{i}}$. 

The lemma is proved.

{\bf Lemma 3.17.}   {\sl  If $G$ is a $Q\sigma T$-group,
 then $G$ is a $\sigma$-$SC$-group. }

 {\bf Proof. }     Suppose that this lemma   is false and let $G$ be a
counterexample of minimal order. Let $D=G^{\frak{N_{\sigma}}}$.
 If $D=1$, then $G$ is $\sigma$-soluble
and so $G$ is a $\sigma$-$SC$-group  by Theorem 2.3. Therefore $D\ne 1$.
Let $R$ be  a minimal normal subgroup of $G$ contained in $D$. Then $G/R$
is an $Q\sigma T$-group by Lemma 3.13. Therefore the choice of $G$
implies that $G/R$ is a  $\sigma$-$SC$-group.    Since
 $(G/R)^{{\mathfrak{N}}_{\sigma }}=D/R$ by Lemmas 2.9 and 2.27,
every
chief  factor of $G/R$   below $D/R$ is simple. Hence
 every
chief  factor of $G$ between  $G^{{\mathfrak{N}}_{\sigma }}$   and $R$ is
simple.
Therefore, in view of the Jordan-H\"{o}lder
theorem for groups with operators \cite[Ch. A, Theorem 3.2]{DH},
it is enough to show that   $R$ is   simple.

 Suppose
that this is  false and let $L$ be
a minimal normal subgroup of $R$.  Then $1  < L  < R$ and $L$
 is $\sigma$-subnormal in $G$,  so $L$ is  modular in $G$ since
$G$ is  an $Q\sigma T$-group.    
  Moreover, $L_{G}=1$ and so every chief factor of $G$ below $L^{G}$ 
 is cyclic by  \cite[Theorem 5.2.5]{Schm}. But
$L^{G}=R$   and so $|R|=p$, a
contradiction. Thus $G$ is a $\sigma$-$SC$-group.

  The lemma is proved.

{\bf Lemma 3.18.}  {\sl If $G$ is a $\sigma$-soluble $Q\sigma T$-group, then 
every Hall $\sigma _{i}$-subgroup $H$ of $G$ is an $M$-group for all $i$.}

{\bf Proof.} By Theorem 3.3,  $G=D\rtimes L$, where $D$   is an abelian  Hall
 subgroup of $G$ of odd order,  $L$ is  $\sigma$-nilpotent and 
  the lattice
of all subgroups ${\cal L}(L)$ of $L$ is modular,  every subgroup of  $D$ is normal
 in $G$, and     $O_{\sigma _{j}}(D)$ has
a normal complement in a Hall $\sigma _{i}$-subgroup of $G$ for all $j$.

It follows that $H_{i}=O_{\sigma _{i}}(D) \times S$, where $O_{\sigma _{i}}(D)$ is a Hall 
$\sigma _{i}$-subgroup of $G$ and $S\leq L^{x}$ for some $x\in G$.  Then  
 $O_{\sigma _{i}}(D)$  and $S\simeq SD/D \leq G/D\simeq L $ are have modular 
   lattices of  the subgroups, so the lattice of all subgroups ${\cal L}(H)$ of $H$ is modular
by \cite[Theorem 2.4.4]{Schm}.

 The lemma is proved.

{\bf Proof of Theorem 3.6.}  First assume that $G$ is a $\sigma$-full
   $Q\sigma T$-groups and let $D=G^{{\mathfrak{S}}_{\sigma}}$  be the
 $\sigma$-soluble residual of $G$. 
 Then  $O^{\sigma _{i}}(D)=D$ for all $i\in I$.  
  Moreover, $G/D$ is a $\sigma$-soluble $Q\sigma T$-group by Lemma 3.13,
 so 
$G/D$ is a $\sigma$-supersoluble by Theorem 3.3 and hence, in fact, 
 $D=G^{{\mathfrak{U}}_{\sigma}}$  is the $\sigma$-supersoluble residual of $G$.

 From Lemma 3.17 it
 follows that $G$ is a $\sigma$-$SC$-group. 
Therefore, if $D\ne 1$,    
  $G$ has a Robinson $\sigma$-complex   $(D, Z(D); U_{1},  \ldots , U_{k})$ 
  by Theorem 2.37 and, in view of Lemma   3.14,  
    for any set   $\{i_{1}, \ldots , i_{r}\}\subseteq \{1, \ldots , k\}$, where
 $1\leq r  < k$, $G$ and $G /U_{i_{1}}'\cdots U_{i_{r}}'$ satisfy
 ${\bf N}_{\sigma _{i}}$ for all $\sigma _{i}\in \sigma (Z(D))$ and 
 ${\bf P}_{\sigma _{i}}$ for all $\sigma _{i}\in \sigma (D)$.
  Therefore the necessity of the condition of the theorem holds.

 Now, assume, arguing 
by contradiction, that $G$ is non-$Q\sigma T$-group of minimal order 
satisfying Conditions (i),  (ii), and (iii). Then  $D\ne 1$.
 We consider a $\sigma$-subnormal
 non-modular subgroup   $U$  of $G$  of minimal order. 
  
First we  show  that 
 $U$ is $\sigma$-supersoluble.  
 Asssume that $k=1$, that is, $D=U_{1}=D'$. Then $Z(D)\leq \Phi (D)$
 and 
$(U\cap D)Z(D)/Z(D)$ is a $\sigma$-subnormal subgroup by Lemma 2.30(1)(2)
  of a non-$\sigma$-primary simple group $D/Z(D)$. 
 Hence $U\cap D\leq Z(D)$, so $U\cap D=U\cap  Z(D)$.
 Therefore 
 every chief factor of $U$  below $U\cap D$ is cyclic. On the other hand, $U/(U\cap D)
\simeq UD/D$ is $\sigma$-supersoluble  by Theorem 3.3 and Proposition 2.28(ii), so $U$ is 
$\sigma$-supersoluble.

Now let $k\ne 1$. We show that the hypothesis holds for $G/U_{i}'$
 for all $i=1, \ldots , k$. We can assume   without loss of generality that $i=1$.
  Let $N=U_{1}'$.
   Then
 $ (G/N)/(D/N)\simeq G/D$ is a $\sigma$-soluble $M\sigma T$-group, so
 Condition (i) holds for $G/D$.  From Lemma
 2.41(1) it follows that   
 $$(D/N, Z(D/N); U_{2}N/N, \ldots , U_{k}N/N)$$ is
 a Robinson $\sigma$-complex of $G/N$ and 
  $$Z(D/N)=U_{1}/N=Z(D)N/N\simeq Z(D)/(Z(D)\cap N).$$
   Moreover, if
 $\{i_{1}, \ldots , i_{r}\}\subseteq \{2, \ldots , k\}$, where $2\leq r  < k$,
 then  the quotients
 $G/N=G/U_{1}'$ and $$(G/N) /(U_{i_{1}}N/N)'\cdots (U_{i_{r}}N/N)'=
(G/N)/(U_{i_{1}}'\cdots U_{i_{r}}'U_{1}'/N)\simeq G/U_{i_{1}}'\cdots U_{i_{r}}'U_{1}'$$
 satisfy
 ${\bf N}_{\sigma _{i}}$ for all
 $\sigma _{i}\in  \sigma (Z(D/N))\subseteq \sigma (Z(D))$ 
  and ${\bf P}_{\sigma _{i}}$ for all
 $\sigma _{i}\in  \sigma (D/N)$ by Condition (iii), so Conditions (ii) and (iii)
 hold for $G/D$.   
Therefore the hypothesis holds for $G/N=G/U_{1}'$, so the 
$\sigma$-subnormal subgroup $UU_{1}'/U_{1}'$ of $G/N$  is modular in
  $G/U_{1}'$ by the choice of $G$.  Finally,  
from Lemma 2.27(2) and Proposition 2.28(2)
 it follows that 
  $$(G/U_{1}')^{\mathfrak{S}_{\sigma}}=U_{1}'G^{\mathfrak{S}_{\sigma}}/U_{1}'
=DU_{1}'/U_{1}'=
 U_{1}'G^{\mathfrak{U}_{\sigma}}/U_{1}'= (G/U_{1}')^{\mathfrak{U}_{\sigma}}.$$
 Hence $U$ is  $\sigma$-supersoluble by
 Lemma 3.15(1).

Next we show that  $U$ is a  $p$-group for some prime $p\in \sigma _{i}$.

Clearly, $(G/N)^{\mathfrak{S}_{\sigma}}= (G/N)^{\mathfrak{U}_{\sigma}}$ 
and so, since $U$ is $\sigma$-supersoluble,
 it is enough  to show that the hypothesis holds
on
$G/N$ for every  nilpotent normal  subgroup $N$ of $G$ by Lemma 3.15(2). 
 It is clear also that $$(G/N)/(DN/N)\simeq G/DN\simeq (G/D)/(DN/D)$$ is
 a $\sigma$-soluble $Q\sigma T$-group, where $|DN/N|\ne 1$.  
 Moreover,  $Z(DN/N)= Z(D)N/N$ and
$$(DN/N,  Z(D)N/N; U_{1}N/N, \ldots ,
U_{k}N/N)$$ is  a Robinson $\sigma$-complex of $G/N$ by Lemma 2.41(2).
In view of Lemmas 2.42 and 3.16,  $G/N$   satisfies ${\bf N}_{\sigma _{i}}$ for all
$$\sigma _{i}\in \sigma (Z(D)N/N)\subseteq \sigma (Z(D))$$ and 
${\bf P}_{\sigma _{i}}$ for all
$$\sigma _{i}\in \sigma (DN/N)\subseteq \sigma (D.)$$

Similarly, for any set 
 $\{j_{1}, \ldots , j_{r}\}\subseteq \{1, \ldots , k\}$, where
 $1\leq r  < k$,  $$(G/N) /(U_{j_{1}}N/N)'\cdots (U_{j_{r}}N/N)'=
(G/N) /(U_{j_{1}}'\cdots U_{j_{r}}'N/N)$$$$\simeq G/U_{j_{1}}'\cdots U_{j_{r}}'N\simeq
   (G/U_{j_{1}}'\cdots U_{j_{r}}')/(U_{j_{1}}'\cdots
U_{j_{r}}'N/U_{j_{1}}'\cdots U_{j_{r}}')
$$  satisfies ${\bf N}_{\sigma _{i}}$ for all
$\sigma _{i}\in \sigma (Z(D)N/N)$ since $$U_{j_{1}}'\cdots
U_{j_{r}}'N/U_{j_{1}}'\cdots U_{j_{r}}'\simeq N/(N\cap U_{j_{1}}'\cdots
U_{j_{r}}')$$ is nilpotent.
 Therefore   the hypothesis holds
on $G/N$, so $U$ is a $p$-group for some prime $p$.

Since $U$ is $\sigma$-subnormal in $G$ by hypothesis,
 for some $p\in \sigma _{i}$ we have $U \leq U^{G}   \leq O_{\sigma _{i}}(G)$ by
 Lemma 2.30(5).

We show that  $U$ modular in $\langle x, U\rangle $ for 
 all elements $x$ of $G$ of prime power order $q^{a}$. 

First assume $q\in \sigma _{i}'$ and we show in this case that 
$x$  induces  a
 power automorphisms in  $O_{\sigma _{i}}(G)$.

Assume that  $O_{\sigma _{i}}(G)\cap D=1$. Since $G/D$ is a $Q\sigma T$-group
by Condition (i), it follows that  $\sigma _{i}'$-elements  of $G/D$ induce
 power automorphisms in  $O_{\sigma _{i}}(G/D)$. Then from the $G$-isomorphism 
$U^{G}\simeq U^{G}D/D$ it follows that  $x$ induces
a  power automorphisms in  $O_{\sigma _{i}}(G)$.

 Now assume that
 $E:=O_{\sigma _{i}}(G)\cap D\ne 1. $  Then $EZ/Z$ is a $\sigma$-soluble normal
 subgroup of $D/Z$, so $E\leq Z$ since $O^{\sigma _{j}}(D/Z)=D/Z$ for all $j$.
Hence $\sigma _{i}\in \sigma (Z)$, so $G$ sasisfies 
 ${\bf N}_{\sigma _{i}}$, which implies   that $x$ 
  induces
 power automorphisms in  $O_{\sigma _{i}}(G)$.  
Therefore  the subgroup $U$ modular in $\langle x, U\rangle $.
                    
Now let $q\in \sigma _{i}$. If $\sigma _{i}\cap \pi (D)\ne \emptyset$, then 
$U$ modular in $\langle x, U\rangle $ since  in this case $G$ sasisfies 
 ${\bf P}_{\sigma _{i}}$ by hypthesis. Finally, assume that  $\sigma _{i}\cap \pi (D)=
 \emptyset$. Then $O_{\sigma _{i}}(G)\cap D=1$.  
Let $E$ be a minimal supplement to $D$ in $G$. Then $E\cap D\leq \Phi(E)$, so $E$ is
 $\sigma$-soluble. Hence $E$ has a Hall $\sigma _{i}$-subgroup $H$. Then  $H$ is 
a Hall $\sigma _{i}$-subgroup of $G$ since $|G:E|=|D/(D\cap E)|$ is a
 $\sigma _{i'}$-number. Therefore $U\leq  O_{\sigma _{i}}(G)\leq H^{y}$, where 
$y$ is such that $x\leq H^{y}$, 
 and $H^{y}\simeq DH^{y}/D$ is a Hall $\sigma _{i}$-subgroup of $G/D$, so $H^{y}$ is an 
$M$-group by Lemma 3.18. But then $U$ modular in $\langle x, U\rangle $. Therefore $U$ 
modular in $\langle x, U\rangle $ for every element $x\in G$ of prime power order.
 It follows that $U$ is modular in $G$ by \cite[Lemma~5.1.13]{Schm}. This contradiction completes the
 proof of the fact that  the sufficiency of the condition of the theorem holds.

The theorem is proved.

The class of all $Q\sigma T$-groups is much wider than the class of
 all $PT$-groups. Indeed, a non-soluble  group $G$
 of order 60 is a $PT$-group and this group
 is not  a $Q\sigma T$-group, where $\sigma =\{\{2, 3, 5\}, \{2, 3, 5\}'\}$,
 since its  subgroups of order 2 are $\sigma$-subnormal but not quasinormal in  $G$.

 Let $\pi \subseteq \Bbb{P}$. Then we say \emph{$G$  satisfies}  
 ${\bf t}_{\pi}$  if whenever $N$ is  a soluble normal
subgroup of $G$, every subgroup of $O_{\pi}(G/N)$ is normal in $G/N$.

By making slight changes to the
 proof of Theorem 3.6 and using this theorem,
 we can prove the following result, which answers Question 3.4(2).

{\bf Theorem 3.19}  (Safonova, Skiba \cite{preprI})
 {\sl Suppose  that $G$ is a $\sigma$-full group.
 A group  $G$ is a  $T_{\sigma}$-group  if 
 and   only if  $G$   satisfies}  
 ${\bf t}_{\sigma _{i}}$   for all   $i\in I$ 
  and every non-$\sigma$-primary chief factor  of $G$ is simple.

{\bf Corollary 3.20}  (Ballester-Bolinces, Beidleman and Heineken    \cite{Ball})
 {\sl  A group  $G$ is a  $T$-group  if 
 and   only if  $G$   satisfies}  
 ${\bf t}_{p}$   for all   primes $p$ ,
  and every non-abelian chief factor  of $G$ is simple.

\section{Groups in which $\sigma$-quasinormality is a transitive relation}

Recall that a subgroup $A$ of $G$ is $\sigma$-quasinormal in $G$ if $A$
 is $\sigma$-subnormal and   modular in $G$.

 We say that   $G$ is an \emph{$M\sigma T$-group} if
 $\sigma$-quasinormality is
a transitive relation on $G$, that is,    if $H$ is a $\sigma$-quasinormal subgroup 
 of $K$ and $K$ is a $\sigma$-quasinormal subgroup 
 of $G$, then $H$ is a $\sigma$-quasinormal subgroup of $G$. 
 
The following open problem
 is one of the most intriguing and difficult problems in the theory of
 $\sigma$-properties  of a group.

{\bf Question 4.1.} 
 {\sl What is the structure of  $M\sigma T$-groups?}

 In this and the following sections,
 we discuss two special cases of this problem.

First note  that in the case when $\sigma =\{\mathbb {P}\}$, this problem 
is the following old open question.

{\bf Question 4.2.} {\sl What is the structure of  \emph{$MT$-groups},
 that are, groups 
 $G$ in which modularity is   a transitive
 relation on $G$, that is,  if $H$ is a modular subgroup of $K$ and $K$ is a modular
 subgroup   of $G$, then $H$ is a modular subgroup of $G$?}

The following special case of Proble. 4.1 is also not explored.

{\bf Question 4.3.}  {\sl What is the structure of
  a $\sigma$-soluble $M\sigma T$-groups?}

The solution of Problem 4.2 in the class of soluble groups was found 
by Frigerio in \cite{A. Frigerio}.

{\bf Theorem  4.4}  (Frigerio \cite{A. Frigerio}, Zimmermann \cite{mod}).
  {\sl A soluble group $G$  is an $MT$-group 
if and only if $G$ is a group with modular latice of all subgroups ${\cal L}(G)$.}

In a later paper \cite{mod},  Zimmermann found a new proof of this result.

In the class of soluble groups,
 the solution of  Problem 4.1 was given for the most general case
 in the recent papers \cite{preprI, matem}.

{\bf Theorem 4.5 }(Zhang, Guo, Safonova, Skiba  \cite{preprI, matem}).  {\sl If  $G$ is a  soluble 
 $M\sigma T$-group 
 and  $D=G^{\frak{N_{\sigma}}}$, then  
    the following conditions hold:} 

(i) {\sl $G=D\rtimes M$, where $D$   is an abelian  Hall
 subgroup of $G$ of odd order, $M$ is a $\sigma$-nilpotent  and the   lattice 
of all subgroups 
${\cal L}(M)$ of $M$  is modular,  }

(ii) {\sl  every element of $G$ induces a
 power automorphism in $D$,   }

(iii) {\sl  $O_{\sigma _{i}}(D)$ has 
a normal complement in a Hall $\sigma _{i}$-subgroup of $G$ for all $i$.}

{\sl Conversely, if  Conditions (i), (ii) and (iii) hold for
  some subgroups $D$ and $M$ of
 $G$, then $G$ is  a soluble $M\sigma T$-group.}

From Theorems  3.5 and 4.5 we get the following interesting result.

{\bf Theorem 4.6 }(Zhang, Guo, Safonova, Skiba  \cite{preprI, matem}).  {\sl Let $G$ be a soluble group. 
Then $G$ is a $W\sigma T$-group  if and only if $G$ is an $M\sigma T$-group}.

In the case $\sigma =\{\mathbb{P}\}$ we get
 from this theorem the following known result.

{\bf Corollary 4.7}  (Frigerio \cite{A. Frigerio}, Zimmermann \cite{mod}).
  {\sl A soluble group $G$  is an $MT$-group 
if and only if $G$ is a group with modular latice of all subgroups ${\cal L}(G)$.}

In the classical case   $\sigma =
\sigma ^{1}=\{\{2\}, \{3\}, \{5\},  \ldots  \}$ we get from Theorem 4.5
 the following   well-known result.

{\bf Corollart 4.8 } (Zacher \cite{zaher}). {\sl  A group
  $G$ is a  soluble  $PT$-group
 if and only if the following conditions are satisfied:}
 
(i) {\sl   the nilpotent residual $L$ of $G$ is an abelian Hall subgroup of odd order,}

(ii) {\sl  $G$ acts by conjugation on $L$ as a group power automorphisms, and }

(iii) {\sl every subgroup of  $G/L$ is quasinormal in $G/D$. }

Nevertheless, the answer to Question 4.3 is  unknown.

The next lemma is a corollary of  general properties  of modular
 subgroups \cite[p. 201]{Schm} and $\sigma$-subnormal subgroups (see Lemma 2.30).

{\bf Lemma 4.9.}  {\sl Let $A$, $B$  and
$N$ be subgroups of $G$, where $A$ is $\sigma$-quasinormal and $N$ is
normal in $G$.}

(1) {\sl The subgroup $A\cap B$ is  $\sigma$-quasinormal  in   $B$.}

(2) {\sl The subgroup  $AN/N$
 is $\sigma$-quasinormal  in  $G/N$}.

(3) {\sl If $N\leq B$ and $B/N$ is  $\sigma$-quasinormal  in  $G/N$, then 
$B$  is
$\sigma$-quasinormal  in  $G$. }

{\bf Lemma 4.10.} {\sl A subgroup $A$ of $G$ is a maximal $\sigma$-quasinormal in $G$
 if and only if either    $G/A=G/A_{G}$ is a 
 simple $\sigma$-primary group  or
 $G/A_{G}$ is a $\sigma$-primary non-abelian group  of order $pq$ for primes
 $p$ and $q$.}

{\bf Proof. }  Assume that  $A$ is a maximal $\sigma$-quasinormal in $G$. In view of
 Theorem 3.3(2), $G/A_{G}$ is a $\sigma _{i}$-group for some $i$, so every subgroup of $G$ 
between $A_{G}$ and $G$ is $\sigma$-subnormal in $G$ by Lemma 2.30(5).
On the other hand, $U/A_{G}$ is modular in $G$ if and only if $U$ is modular in $G$
 by \cite[Page 201, Property (4)]{Schm}. Therefore, in fact, $A$ is a maximal modular 
subgroup of $G$. Hence  either $A\trianglelefteq G$ and    $G/A=G/A_{G}$ is a 
 simple $\sigma _{i}$-group  or
 $G/A_{G}$ is a non-abelian group   of order $pq$ for primes $p, q \in \sigma _{i}$ 
by \cite[Lemma 5.1.2]{Schm}.

Finally, if  $G/A_{G}$ is a $\sigma$-primary non-abelian group   of order $pq$ for primes
 $p$ and $q$, then $A$ is a maximal subgroup of $G$ and $A/A_{G}$
 is, clearly,  modular in $G/A_{G}$, so $A$ is a maximal modular subgroup of $G$
 by \cite[Page 201, Property (4)]{Schm},  so $A$ is 
 a maximal $\sigma$-quasinormal of $G$.
 Similarly,
 if $G/A=G/A_{G}$ is a 
 simple non-abelian $\sigma$-primary group, then $A$ is a maximal modular
 subgroup of $G$ by \cite[Lemma 5.1.2]{Schm}, so $A$ is 
 a maximal $\sigma$-quasinormal of $G$.

The lemma is proved.

We say that a  subgroup $A$ of $G$ is said to be  \emph{$\sigma$-subquasinormal} in $G$ 
if     
   there is a subgroup chain  $A=A_{0} \leq A_{1} \leq \cdots \leq
A_{n}=G$  such that  $A_{i-1}$ is $\sigma$-quasinormal in $ A_{i}$ 
 for all $i=1, \ldots , n$.

It is clear that $G$ is an $M\sigma T$-group if and only if every its 
$\sigma$-subquasinormal subgroup is $\sigma$-quasinormal in $G$.

We use $\frak{A}^{*}$ to denote the class of all abelian groups of squarefree exponent;
$G^{{\frak{A}^{*}}}$ is the intersection of all  normal subgroups $N$ of $G$ 
with $G/N\in {\frak{A}^{*}}$.  It is clear that $\frak{A}^{*}$ is a hereditary formation, 
so $G/G^{{\frak{A}^{*}}}\in \frak{A}^{*}$.

The following lemma is a corollary of Lemmas 1 and 4 in  \cite{mod} and Lemma
 2.30.

{\bf Lemma 4.11.} {\sl Let  $A$,  $B$ and $N$ be subgroups of $G$,
 where 
  $A$  is  $\sigma$-subquasinormal  $G$ and $N$ is normal in $G$.  }

(1) {\sl $A\cap B$    is   $\sigma$-subquasinormal  in   $B$}.

(2) {\sl $AN/N$ is
 $\sigma$-subquasinormal  in $G/N$. }

(3) {\sl If $N\leq K$ and $K/N$ is   $\sigma$-subquasinormal 
in $G/N$, then $K$ is   $\sigma$-subquasinormal  in $G.$}

 (4)    {\sl   $A^{{\frak{A}^{*}}}$    is subnormal in $G$.}

{\bf Lemma 4.12} (See \cite[Lemma 5.1.13]{Schm}).  {\sl  A subgroup $M$ of $G$ is modular 
in $G$ if and only if $M$ is modular in $\langle x, M\rangle$ for every 
 element of prime power order  $x\in G$.}

{\bf Proposition 4.13.} {\sl  Suppose that  $G$ is a  soluble $PT$-group and
 let $p$ be a prime.
 If every submodular $p$-subgroup of $G$ is modular in $G$, then every
 $p$-subgroup of $G$ is modular in $G$. In particular, if
 every submodular subgroup of a soluble 
$PT$-group $G$ is modular in $G$, then $G$  is an $M$-group.}

{\bf Proof.}  Assume that this proposition is false and let $G$ be a
 counterexample of minimal order. Then, by \cite[Theorem 2.1.11]{prod},
  the following conditions are satisfied:
the nilpotent residual $D$ of $G$ is a Hall subgroup,
 $G$ acts by conjugation on $D$ as a group power automorphisms, and
  every subgroup of  $G/D$ is quasinormal in $G/D$. In particular, $G$ is supersoluble.

 Let $M$ be a complement to $D$ in $G$ and  $U$ a
 non-modular $p$-subgroup of $G$ of minimal order. Then $U$ is not submodular in $G$ and
 every maximal subgroup of $U$ 
is modular in $G$, so $U$ is a cyclic group by  \cite[p. 201, Property (5)]{Schm}.
Let $V$ be the maximal subgroup of $U$.  Then $V\ne 1$ since every subgroup of
prime order of  a supersoluble group is submodular by  \cite[Lemma 6]{mod}. 
We can assume without loss of generality that $U\leq M$ since $M$ is a
 Hall subgroup of $G$.

(1) {\sl  If $R$ is a normal $p$-subgroup of $G$, then
 every $p$-subgroup of $G$ containg $R$ is modular in $G$. In particular,
 $U_{G}=1$ and so  $U\cap D=1.$}

Let $L/R$ be a submodular $p$-subgroup of $G/R$.
 Then $L$ 
is a submodular $p$-subgroup of  $G$ by  Lemma 4.9(3), so $L$
 is modular in $G$ by hypothesis. Hence $L/R$ is modular in $G/R$ by 
\cite[p. 201, Property (4)]{Schm}.  Hence  the hypothesis holds for $G/R$.
Therefore every $p$-subgroup $S/R$ of $G/R$ is modular in $G/R$ by the choice of $G$, so 
$S$ is modular in $G$ by \cite[p. 201, Property (4)]{Schm}.

(2) {\sl If $K$ is a proper submodular subgroup of $G$, then every $p$-subgroup $L$
 of $K$ is modular  in $G$, so every proper subgroup of $G$ contaning $U$ is not 
submodular in $G$. } 

First note that $K$ is a $PT$-group by \cite[Corollary~2.1.13]{prod} and if $S$ is a
 submodular $p$-subgroup of $K$,
then $S$ is submodular in $G$ and so $S$ is modular in $G$.
 Hence $S$ is modular in $K$. 
 Therefore the hypothesis holds for $K$, so  every $p$-subgroup $L$
 of $K$ is modular  in $K$ by the choice of $G$. Hence 
 so $L$ is modular in $G$ by hypothesis.

(3) $DU=G$ (This follows from Claim (2) and the fact that every subgroup of $G$
 containing $D$ is subnormal in $G$).      

(4) {\sl  $V$ is not subnormal in $G$.}

Assume that $V$ is subnormal in $G$. Then $V$ is quasinormal in $G$ by
 Theorem A since $V$ is
 modular  in $G$.  
 Therefore $1 < V \leq R=O_{p}(Z_{\infty}(G))$ by \cite[Corollary 1.5.6]{prod} 
 since $V_{G}=1=U_{G}$ by Claim (1). But $R\leq U$ since   $U$ is a Sylow $p$-subgroup
 of $G$  by Claim (3), hence  $R=V=1$ and
 so $|U|=p$, a contradiction.   Therefore we have (4).

(5) {\sl $G=V^{G}\times K$, where $V^{G}$ is a non-abelian $P$-group of order prime
 to $|K|$}  (Since $V_{G}=1$, this follows from Claim (4) and 
 Lemma~2.2.

{\sl Final contradiction.} From Claim (5) it follows that  $U\leq V^{G}$, so $U$ is 
submodular in $G$ by \cite[Theorem 2.4.4]{Schm}). This final contradiction
 completes the proof of the result. 

The proposition is proved.

{\bf Proof of Theorem 4.4.}   Let  $\sigma (G)=\{\sigma _{1}, \ldots, \sigma _{t}\}$.
  
 First suppose that $G$ is  a soluble $M\sigma 
T$-group.  Then  $G$ has a Hall $\sigma _{i}$-subgroup  $H_{i}$ for all
 $i=1, \ldots, t$.

We show that Conditions (i), (ii), and  (iii)  hold for $G$.  Assume that
this is false and let $G$ be a counterexample of minimal order.

(1) {\sl If $R$ is a non-identity normal subgroup of $G$, then
 Conditions (i), (ii), and  (iii) hold for  $G/R$.}

If $H/R$ is a $\sigma$-subquasinormal subgroup of $G/R$, then $H$ is
 $\sigma$-subquasinormal
 in $G$ by Lemma 4.11(3), so $H$ is $\sigma$-quasinormal in $G$ by hypothesis
 and hence $H/R$ is a $\sigma$-quasinormal  in $G/R$ by Lemma 4.9(2). Therefore
 $G/R$ is an 
$M\sigma T$-group, so we have (1) by the choice of $G$.

(2) {\sl If $E$ is a proper $\sigma$-subquasinormal subgroup of $G$, then 
 $E^{\frak{N_{\sigma}}}\leq D$ and Conditions (i) and (ii)  hold for  $E$.   }

Every  $\sigma$-subquasinormal   subgroup $H$ of $E$ is $\sigma$-subquasinormal 
in $G$,  so $H$ is  $\sigma$-quasinormal 
in $G$ by hypothesis.    Therefore the hypothesis holds for $E$, 
so Conditions (i) and  (ii)  hold for  $E$ by the choice of $G$.  
Moreover,  since $G/D\in {\frak{N_{\sigma}}}$ 
  and ${\frak{N_{\sigma}}}$ is a hereditary  class by Lemma 2.9, 
$E/E\cap D\simeq ED/D\in {\frak{N_{\sigma}}}$  and so    
$E/E\cap D    \in {\frak{N_{\sigma}}}$. Hence $$E^{\frak{N_{\sigma}}}\leq 
E\cap D\leq D.$$ 

(3) {\sl $D$ is nilpotent and $G$ is supersoluble.}

 Let $R$ be a minimal   normal subgroup of $G$. Then 
$RD/R=(G/R)^{{\frak{N_{\sigma}}}}$ is abelian by  Lemma 2.27 and Claim (1).
 Therefore $R\leq D$, $R$   
is the unique minimal   normal subgroup of $G$ and $R\nleq \Phi (G)$. Hence  
   $R=C_{G}(R)=O_{p}(G)=F(G)$ for some  
$p\in \sigma _{i}$ by \cite[Ch. A, 13.8(b) and 15.2]{DH}

  Let $V$ be a maximal 
subgroup of $R$. Then $V_{G}=1$ and $V$ is $\sigma$-subquasinormal 
in $G$,  so $H$ is  $\sigma$-quasinormal in $G$. Hence $R=V^{G}$ is a group of order $p$ 
by \cite[Theorem  5.2.3]{Schm}, so $ G/R=C_{G}(R)$
 is cyclic and  hence $G$ is supersoluble. But 
then $D=G^{{\frak{N_{\sigma}}}}\leq G'\leq F(G)$ and  so $D$   is nilpotent.

(4) {\sl If $H/K$ is   $\sigma$-nilpotent, where $K\leq H$ are normal subgroups of $G$,
then $H/K$ is an $M$-group.}

Let $U/K$ be any submodular subgroup of $H/K$, then  $U/K$ is submodular in $G/K$ and so 
$U$ is  is submodular in $G$. 
 On the other hand,
$U/K$ is  $\sigma$-subnormal in $G/D$ by Lemma 2.10,
 so $U$ is $\sigma$-subnormal  in $G$ by 
Lemma 2.30(5). Therefore $U$ is $\sigma$-subquasinormal in $G$
 and so, by hypothesis,  
$U$ is $\sigma$-quasinormal in $G$. It follows that  $U/K$ is
 $\sigma$-quasinormal in $G/K$ and so in $H/K$ by Lemma 4.9(1)(2), so every 
 submodular subgroup of $H/K$
 is modular $H/K$. Therefore $H/K$ is an $M$-group by Proposition 4.13.

(5) {\sl  $D$ is a Hall subgroup of $G$. }

 Suppose
that this is false and let $P$ be a  Sylow $p$-subgroup of $D$ such
that $1 < P < G_{p}$, where $G_{p}\in \text{Syl}_{p}(G)$.  We can assume 
without loss of generality that $G_{p}\leq H_{1}$.

(a)  {\sl    $D=P$ is  a minimal normal subgroup of $G$. }

Let $R$ be a minimal normal subgroup of $G$ contained in $D$. 
Then   $R$ is a $q$-group    for some prime   
$q$. Moreover, 
$D/R=(G/R)^{\mathfrak{N}_{\sigma}}$  is a Hall subgroup of $G/R$ by
Claim (1).  Suppose that  $PR/R \ne 1$. Then
  $PR/R \in \text{Syl}_{p}(G/R)$. 
If $q\ne p$, then    $P \in \text{Syl}_{p}(G)$. This contradicts the fact 
that $P < G_{p}$.  Hence $q=p$, so $R\leq P$ and therefore $P/R \in 
\text{Syl}_{p}(G/R)$ and we again get  that 
$P \in \text{Syl}_{p}(G)$. This contradiction shows that  $PR/R=1$, which implies that 
  $R=P$  is the unique minimal normal subgroup of $G$ contained in $D$.
 Since $D$ is nilpotent,
 a $p'$-complement $E$ of $D$ is characteristic in 
$D$ and so it is normal in $G$. Hence $E=1$, which implies that $R=D=P$.

(b) {\sl $D\nleq \Phi (G)$.    Hence for some maximal subgroup
 $M$ of $G$ we have $G=D\rtimes M$  }  (This follows  from  Lemma 2.9 since $G$
 is not $\sigma$-nilpotent).

(c) {\sl If $G$ has a minimal normal subgroup $L\ne D$,
 then $G_{p}=D\times (L\cap G_{p})$.
  Hence $O_{p'}(G)=1$. }

Indeed, $DL/L\simeq D$ is a Hall 
subgroup of $G/L$ by Claims (1)  and (a). Hence
  $G_{p}L/L=RL/L$, so $G_{p}=D\times (L\cap G_{p})$.
 Thus  $O_{p'}(G)=1$ since $D < G_{p}$ by Claim (a).

(d)  {\sl   $V=C_{G}(D)\cap M$ is a  normal subgroup of $G$ and 
 $C_{G}(D)=D\times V \leq H_{1}$.  }

In view of  Claims  (a) and 
 (b),  $C_{G}(D)=D\times V$, where $V=C_{G}(D)\cap M$ 
is a normal  subgroup of $G$. Moreover,  $V\simeq DV/D$ is $\sigma $-nilpotent by
 Lemma 2.9.  Let $W$ be a $\sigma 
_{1}$-complement of $V$. Then $W$  is characteristic in $V$ and so it is normal 
in $G$.    Therefore we have  (d) by Claim (c).

(e)  $G_{p}\ne H_{1}$.

Assume that $G_{p}=H_{1}$.  Let $Z$ be a subgroup of order $p$ in $Z(G_{p})\cap D$.
Then $Z$ is $\sigma$-subquasinormal in $G$ by Claim (3), so $Z$ is
 $\sigma$-quasinormal in $G$ and hence  $ O^{\sigma _{1}}(G)=O^{p}(G)\leq N_{G}(Z)$
by Theorem 3.3(v).
Therefore $G=G_{p}O^{p}(G)\leq N_{G}(Z)$, hence $D=Z < G_{p}$ If follows that
 $D < C_{G}(D)$. 
   Then  $V=C_{G}(D)\cap M\ne 1$ is a normal subgroup of $G$ and 
  $V\leq H_{1}=G_{p}$ by Claim (d). Let $L$ be a   minimal
 normal subgroup of $G$ contained in $V$. Then  $G_{p}=D\times L$ is a normal  
elementary abelian subgroup of $G$ of order $p^{2}$ by Claim  (c) and  
 every subgroup of $G_{p}$ is 
normal in $G$ by Theorem 3.3(v). 
Let $D=\langle a \rangle$,  $L=\langle b \rangle$ and $N=\langle ab \rangle$.  
Then $N\nleq D$, so in view of the $G$-isomorphisms
 $$DN/D\simeq N\simeq NL/L= G_{p}/L=DL/L\simeq D $$  we get that 
$G/C_{G}(D)=G/C_{G}(N)$ is a $p$-group since $G/D$ is $\sigma$-nilpotent by Lemma 
2.9.
But then Claim (d) implies that  $G$ is a $p$-group. This 
contradiction shows that we have (e).

{\sl Final contradiction for (5).} In view of Theorem A in \cite{2}, $G$ has a $\sigma 
_{1}$-complement $E$ such that $EG_{p}=G_{p}E$. 

Let $V=(EG_{p})^{{\frak{N}}_{\sigma}}$.  By 
Claim (e), $EG_{p}\ne G$.     On the other hand, since $   D\leq 
EG_{p}$ by Claim (a),  $EG_{p}$ is $\sigma$-subquasinormal in $G$ by Claim (4) and
 Lemma 4.11(3). 
 Therefore   Claim (2) implies that    $V$  is a Hall subgroup of $EG_{p}$ 
 and  $V\leq D$, 
 so  for a Sylow 
$p$-subgroup $V_{p}$ of $V$ we have $|V_{p}|\leq |P| < |G_{p}|$. 
Hence    $V=1$.  
Therefore $EG_{p}=E\times G_{p}$ is $\sigma$-nilpotent and so
 $E\leq C_{G}(D)\leq H_{1}$. Hence $E=1$ and  so $ D =1$, a contradiction.  Thus,   
$D$ is a Hall subgroup of $G$.

(6)  {\sl $H_{i}=O_{\sigma _{i}}(D)\times S$ 
 for each $\sigma _{i} \in \sigma (D) $.}

Since  $D$ is  a nilpotent Hall subgroup of $G$ by Claims (3) and (5),
 $D=L\times N$, where $L=O_{\sigma _{i}}(D)$  
  and $N=O^{\sigma _{i}}(D)$ are Hall subgroups of $G$.  
 First assume that $N\ne 1$.  Then $$O_{\sigma 
_{i}}((G/N)^{\mathfrak{N}_{\sigma}})=O_{\sigma 
_{i}}(D/N)=LN/N$$ has a normal complement $V/N$ in $H_{i}N/N\simeq H_{i}$
  by Claim (1). On the other hand,  $N$ has a complement $S$ in $V$ by 
the Schur-Zassenhaus theorem.   Hence  $H_{i}=.H_{i} \cap LSN=LS$ and $L\cap 
S=1$ since $$(L\cap 
S)N/N\leq (LN/N)\cap (V/N)=(LN/N)\cap (SN/N)=1$$ 
 It is clear that $V/N$ is a Hall subgroup of $H_{i}N/N$, so $V/N $  is 
characteristic in $H_{i}N/N$. On the other hand,  $H_{i}N/N$ is 
normal  in $G/N$ by Lemma 2.9 since $D/N\leq H_{i}N/N$. 
 Hence $V/N $  is  normal in $G/N$.
  Thus $H_{i}\cap V =H_{i}\cap NS=S(H_{i}\cap N)=S
 $ is normal in $H_{i}$, so $H_{i}=O_{\sigma _{i}}(D)\times S$.

Now assume  that  $D=O_{\sigma _{i}}(D)$. Then $H_{i}$ is normal in 
$G$, so   $H_{i}$ is an $M$-group by Claim (4).
Therefore every subgroup $U$ of $H_{i}$ is   $\sigma$-quasinormal and so 
 $\sigma$-normal in $G$ by Theorem B(v). Since $D$ is a normal Hall subgroup of $H_{i}$,
 it has a 
complement $S$ in $H_{i}$ and  $D\leq   O^{\sigma _{i}}(G)\leq N_{G}(S)$ since $S$ is 
$\sigma$-normal in $G$. Hence 
$H_{i}=D\times S=O_{\sigma _{i}}(D)\times S$.

(7) {\sl Every subgroup $H$ of $D$ is normal in $G$. Hence every element of
 $G$ induces a power automorphism in $D$. }

Since $D$ is  nilpotent by Claim (3), it is enough to consider 
the case when $H\leq O_{\sigma _{i}}(D)=H_{i}\cap D$ for some $\sigma _{i}\in \sigma (D)$.
Claim (6) implies that $H_{i}=O_{\sigma _{i}}(D)\times S$. 
It is clear that $H$ is $\sigma$-subquasinormal in $G$, so $H$ is $\sigma$-quasinormal
 in 
$G$. Therefore $H$ is $\sigma$-normal in $G$  by Theorem B(v),
 so   $$G=H_{i}O^{\sigma _{i}}(G)=
(O_{\sigma _{i}}(D)\times S)O^{\sigma _{i}}(G)=SO^{\sigma _{i}}(G)\leq 
N_{G}(H)$$.

 (8) {\sl  If  $p$ is a  prime such that $(p-1, |G|)=1$, then  $p$
does not divide $|D|$. Hence the smallest prime in $\pi (G)$ belongs to $\pi (|G:D|)$.
 In particular,  $|D|$ is odd. }

Assume that this is false.
 Then, by Claim (7),  $D$ has a maximal subgroup $E$ such that
$|D:E|=p$ and  $E$ is normal in $G$. It follows that  $C_{G}(D/E)=G$ since $(p-1, 
|G|)=1$.   
Since
$D$ is a Hall subgroup of $G$ by Claim (5), it has a complement $M$ in $G$. 
Hence 
$G/E=(D/E)\times (ME/E)$, where
$ME/E\simeq M\simeq G/D$ is $\sigma$-nilpotent. Therefore $G/E$ is
$\sigma$-nilpotent by Lemma 2.9. But then $D\leq E$, a contradiction. 
Hence we have (8).

(9)  {\sl  $D$ is abelian.}

In view of Claim 
(7), $D$ is a Dedekind group.  Hence $D$ is abelian since $|D|$ is  odd  by Claim (8).  

 From Claims (4)--(9) we get that Conditions (i), (ii),  and (iii), where $M$ is an $M$-group
  hold  for $G$.

Now, we show that if Conditions (i), (ii),and  (iii), where $M$ is an $M$-group,
  hold for  $G$,
 then
 $G$ is an $M\sigma T$-group.  Assume that
this   is false and let $G$ be a counterexample of minimal order.  
  Then  $D\ne 1$ and, by Lemma 4.12,   
 for some  $\sigma$-subquasinormal 
 subgroup $A$ of $G$ and for
some element $x\in G$ of prime power order $p^{a}$  the subgroup $A$ is not 
modular in $\langle A, x  \rangle $. Moreover, we can assume that every proper 
$\sigma$-subquasinormal subgroup of $A$ is $\sigma$-quasinormal  in $G$.

(*) {\sl   If  $N$ is a minimal normal subgroup of $G$, then
 $G/N$ is an $M\sigma T$-group. } (Since the hypothesis holds for $G/N$, this follows from
 the choice of $G$).

(**) {\sl If $N$ is a minimal normal subgroup of $G$, then $AN$ is
$\sigma$-quasinormal   in  $G$. In particular, $A_{G}=1$.}

Claim (*) implies that $G/N$ is an $M\sigma T$-group. On the other hand,
  by Lemma 4.11(2), 
    $AN/N$ is a $\sigma$-subquasinormal subgroup of 
$G/N$, so  $AN/N$ is a $\sigma$-quasinormal in $G$. Hence we have (**) by 
Lemma 4.11(3). 

(***) {\sl $A$ is a $\sigma _{i}$-group for some $i$.}

From Condition (i) and Claim (**) it follows that $A\cap D=1$,
 so $AD/D\simeq A$ is $\sigma$-nilpotent.
Then for $\sigma _{i}\in \sigma (A)$ we have
 $A=O_{\sigma _{i}}(A)\times O_{\sigma _{i}'}(A)$. Assume that
 $O_{\sigma _{i}'}(A)\ne 1$. Then  $O_{\sigma _{i}}(A)$ and
 $ O_{\sigma _{i}'}(A)$
are $\sigma$-quasinormal  in $G$, so $A$ is  $\sigma$-quasinormal 
 in $G$ by Lemma  2.30(3) and  \cite[Page 201, Property (5)]{Schm}.  This
 contradiction shows that $A=O_{\sigma _{i}}(A)$ is   a $\sigma _{i}$-group.

 {\sl Final contradiction for the sufficiency.}  Since $A$ is $\sigma$-subnormal in $G$, 
from Claim (***)  it follows that   $A\leq H_{i}^{y}$  for all $y\in G$ 
 by Lemma 2.30(5).

From Conditions (i) and (ii) it follows that $H_{i}=(H_{i}\cap D)\times S$
  for some Hall subgroup $S$ of $H_{i}$. Then from 
 $A\leq H_{i}^{y}=(H_{i}^{y}\cap D)\times S^{y}$ it follows that  $A\leq S^{y}$
 for all $y\in G$.  
  Moreover, $S\simeq  SD/D$ by Claim (**) and $H_{i}\cap D $ are
  $M$-groups by Condition (i).  Hence $H_{i}$ is an $M$-group by   
 \cite[Theorem 2.4.4]{Schm}.
 Therefore  $x\not \in H_{i}^{y}$  all $y\in G$, so
 $p\in \sigma_{j}$ for some $j\ne i$.

Let $U= \langle x\rangle$. First assume that $x\in D$, then
 $U\trianglelefteq G$. On the other hand, $A$ is $\sigma$-subnormal in $UA$,
 so $UA=U\times A$ by Lemma 2.30(5). Hence $A$ is modular in
 $\langle x, A\rangle =UA$, a contradiction. Therefore, $x\not \in D$. 

Since $D$ is a Hall subgroup of $G$ and $x\not \in D$, $x\in M^{z}$
 for some $z\in G$. It is
also clear that
$ A\leq S^{y} \leq M^{z}$ for some $y\in G$, where $M^{z}$ is an
 $M$-group by Condition (i), and then   $A$ is  
modular in $\langle A, x  \rangle $.

This final contradiction  
 completes the proof of the fact that $G$ is an $M\sigma T$-group.

                               The theorem is proved.

 \section{Groups in which modularity is transitive}

Recall that a group $G$ is said to be an \emph{$M$-group}  \cite[p. 54]{Schm}
 if the  latice ${\cal L}(G)$, of all subgroups of $G$, is modular.

{\bf Definition 5.1.} We say that   $G$ is
 an \emph{$M T$-group}
 if modularity is a transitive relation on $G$, that is, if $H$ is a modular subgroup 
 of $K$ and $K$ is a modular subgroup 
 of $G$, then $H$ is a modular subgroup of $G$.

The next problem goes back to the paper by Frigerio \cite{A. Frigerio}.

{\bf Problem 5.1.} {\sl What is the structure of an   $MT$-group $G$, that is, a group
 in which  modularity is a transitive relation on $G$?}

 Frigerio proved the following theorem which gives a complete answer to
 this  problem  for the soluble case. 
 
{\bf Theorem 5.3}  (Frigerio \cite{A. Frigerio}).  {\sl A soluble group is an $MT$-group 
if and only if $G$ is a group with modular latice of all subgroups ${\cal L}(G)$.}

A new proof of Theorem 5.3 was obtained in the paper \cite{mod}.

 {\bf Remark 5.4.}  Problem 5.2 is a special case
 of general Propblem 4.2, where $\sigma =\{\mathbb{P}\}$, since in this case every
 subgroup of every group is $\sigma$-subnormal.

 Before continuing, we give a few definitions.

A  subgroup $A$ of $G$  is said to be  \emph{submodular} in $G$ if   
   there is a subgroup chain  $A=A_{0} \leq A_{1} \leq \cdots \leq
A_{n}=G$  such that   $A_{i-1}$ is a modular subgroup of $ A_{i}$ 
  for all $i=1, \ldots , n$.  Thus a group $G$ is an $M T$-group if and only if every
 submodular
 subgroup of $G$ is modular in $G$.

{\bf Remark 5.5.}  It is clear that every subnormal subgroup is
 submodular. On the other hand,
 in view of the above mentioned Ore's result in \cite{5} and Theorem 3.1,
 $G$ is a $PT$-group if and only if 
every its subnormal subgroup is modular.   Therefore every  $MT$-group is a
 $PT$-group.

 In view of Remark 5.5, the following well-known 
 result partially  describes the structure of  insoluble $MT$-groups.

{\bf Theorem 5.6}  (Robinson  \cite{217}).  {\sl $G$ is a $PT$-group if 
 and   only if  $G$  has a normal perfect subgroup $D$ such that:}

(i) {\sl  $G/D$ is a soluble $PT$-group, and }

(i) {\sl if $D\ne 1$, $G$ has a Robinson complex
 $(D, Z(D); U_{1},  \ldots , U_{k})$ and }

(iii) {\sl   for any set  $$\{i_{1}, \ldots , i_{r}\}\subseteq \{1, \ldots , k\},$$ where
 $1\leq r  < k$,  $G$ and $G /U_{i_{1}}'\cdots U_{i_{r}}'$ satisfy
 ${\bf N}_{p}$ for all $p\in \pi (Z(D))$ and
 ${\bf P}_{p}$ 
for all $p \in \pi (D)$. }

Now, recall that $G$ is a non-abelian $P$-group (see \cite[p. 49]{Schm}) if 
 $G=A\rtimes \langle t \rangle$, where $A$ is    an elementary abelian
$p$-group and an element $t$ of
 prime order $q\ne p$ induces a non-trivial power
 automorphism  on $A$. In this case we say that $G$ is a \emph{$P$-group
 of type  $(p, q)$}.

{\bf Definition 5.7.}  We say that  $G$  \emph{satisfies 
 ${\bf M }_{P}$} (\emph{${\bf M }_{p, q}$}, respectively) if  
 whenever $N$ is  a soluble normal
subgroup of $G$ and $P/N$ is a normal non-abelian  $P$-subgroup  (a normal $P$-group
 of type  $(p, q)$, respectively)
 of $G/N$,
 every non-subnormal subgroup of $P/N$ is modular in  $G/N$.

The following theorem answers to Problem 5.1.

{\bf Theorem 5.8}  (Liu,   Guo,   Safonova and Skiba  \cite{preprI, pure}).
  {\sl A group $G$ is an $MT$-group if   
 and   only if  $G$  has a perfect normal subgroup $D$ such that:}

(i) {\sl  $G/D$ is an $M$-group,  }

(ii) {\sl if  $D\ne 1$,  $G$ has a Robinson complex
 $(D, Z(D); U_{1},  \ldots , U_{k})$ and }

(iii) {\sl   for any set $$\{i_{1}, \ldots , i_{r}\}\subseteq \{1, \ldots , k\},$$ where
 $1\leq r  < k$,  $G$ and $G /U_{i_{1}}'\cdots U_{i_{r}}'$ satisfy
 ${\bf N}_{p}$ for all $p\in \pi (Z(D))$,
 ${\bf P}_{p}$ 
for all $p\in \pi (D)$, and ${\bf M}_{p, q}$ for all pears
 $\{p, q\}\cap \pi (D)\ne \emptyset.$}

The following example
 shows that, in general, a $PT$-group may not be an $MT$-group.

{\bf Exapmle 5.9.} (i) Let $\alpha: Z(SL(2, 5))\to Z(SL(2, 7))$ be an isomorphism and let 
$$D:= SL(2, 5) \Ydown SL(2, 7)=(SL(2, 5)\times SL(2, 7))/V,$$
 where $$V=\{(a, (a^{\alpha})^{-1})\mid a\in Z(SL(2, 5))\},$$
  is the direct product  of the groups $SL(2, 5)$ and $SL(2, 7)$ with joint center
 (see \cite[p. 49]{hupp}).    
Let   $M=(C_{7}\rtimes C_{3}) \Yup (C_{13}\rtimes C_{3}$) be
 the direct product  of the groups $C_{7}\rtimes C_{3}$ and $C_{13}\rtimes C_{3}$
 with joint 
factorgroup $C_{3}$  (see \cite[p. 50]{hupp}),
 where 
$C_{7}\rtimes C_{3}$ is  a non-abelian group of
 order 21   and  $C_{13}\rtimes C_{3}$ is a non-abelian group of
 order 39. Finally,  let $G=D\times M$. 
 We show that $G$ satisfies the conditions in Theorem 5.6. 

 It is clear that $D=G^{\mathfrak{S}} $ is the soluble
 residual  of $G$ and $M\simeq G/D$ is a   soluble $PT$-group.  
 In view of \cite[I, Satz 9.10]{hupp}, $D=U_{1}U_{2}$ and 
$U_{1}\cap U_{2}=Z(D)=\Phi (D)$, where  $U_{i}$ is normal in $D$, 
 $U_{1}/Z(D)$ is a  simple group of order 60, and 
$U_{2}/Z(D)$ is a  simple group of order 168. Hence $(D, Z(D); U_{1}, U_{2})$ is
 a Robinson    complex   of $G$ and
 the subgroup $ Z(D)$ has order 2 and  $Z(D)\leq Z(G)$.
  Therefore Conditions (i) and (ii) hold for $G$. It is not difficult to show 
 that for every prime 
$r$ dividing $|G|$ and for $O_{r}(G/N)$, where $N$ is a normal
 soluble subgroup of $G$, we have  $|O_{r}(G/N)|\in \{1, r\}$, so Condition (iii)
 also holds for $G$.  
 Therefore  $G$ is a $PT$-group by Theorem~5.6.

Now we show that  $G$ is not an $MT$-group.  First note that $M$ has
 a subgroup $T\simeq C_{7}\rtimes C_{3}$ and $|M:T|=13$. Then $T$ is a maximal
 subgroup of $M$  and  $M/T_{M}\simeq C_{7}\rtimes C_{3}$.
 Hence a subgroup $L$ of $T$ of order 3 is modular in $T$ and $T$ is modular in $M$ by
\cite[Lemma 5.1.2]{Schm},   so $L$  is submodular in  $G$. Finally,
  $L$ is not modular in $M$ by Lemma 5.11 below.
 Therefore $G$ is not an $MT$-group by Theorem~5.8.

(ii)  The group $D \times (C_{7}\rtimes C_{3})$  is an $MT$-group by Theorem~5.8.

{\bf Lemma 5.10.} {\sl Let  $A$,  $B$ and $N$ be subgroups of $G$, where 
  $A$
is submodular  and $N$ is normal in $G$.  Then:}

(1) {\sl $A\cap B$    is  submodular in   $B$},

(2) {\sl $AN/N$ is
submodular in $G/N$, }

(3) {\sl if $N\leq K$ and $K/N$ is  submodular
in $G/N$, then $K$ is  submodular in $G,$}

 (4)    {\sl  $A^{{\frak{A}^{*}}}$    is subnormal in $G$, }

 (5)    {\sl  if $G=U_{1}\times \cdots \times U_{k}$, where $U_{i}$ is a simple
 non-abelian group, then $A$ is normal in $G$. }

{\bf Proof.}  Statements (1)--(4) are proved in \cite{mod}.

(5)  Let $E=U_{i}A$ and $A\ne 1$. Then $A$ is submodular in $E$ by Part (1), so
 there is a subgroup chain 
$$A=E_{0} < E_{1} < \cdots < E_{t-1} < E_{t}=E$$ such that $E_{i-1}$  is a maximal
modular subgroup of  $E_{i}$ for all $i=1, \ldots, t$ and for $M=E_{t-1}$ we have 
$M=A(M\cap U_{i})$ and, by  \cite[Lemma 5.1.2]{Schm}, either
 $M=E_{t-1}$ is  a maximal normal subgroup of $E$ or $M$
 is a maximal subgroup of $E$ such  that 
$E/M_{E}$ is a non-abelian group of order $qr$ for primes  $q$ and $r$.
In the former case we have $M\cap U_{i}=1$, so $A=M$ is normal in $E$. The second case
 is impossible 
since $E$ has no a quotient of order $qr$. Therefore  $U_{i}\leq N_{G}(A)$ for all $i$,
 so 
$G\leq N_{G}(A)$. Hence we have (5).

The lemma is proved.

{\bf Lemma  5.11} (See Lemma 5.1.9 \cite{Schm}). {\sl Let $M$ be a modular subgroup of $G$ 
of prime power order. If $M$ is not quasinormal in $G$, then 
$$G/M_{G}=M^{G}/M_{G}\times K/M_{G},$$ where $M^{G}/M_{G}$ is a non-abelian
 $P$-group of order prime to $|K/M_{G}|$.   }

Recall that a group $G$ is said to be an \emph{$SC$-group}
 if every chief factor of $G$ is simple \cite{217}.

{\bf Lemma 5.12.}  {\sl Let $G$ be a  non-soluble $SC$-group  and suppose that $G$ has a 
  Robinson complex
 $$(D, Z(D); U_{1}, \ldots ,
U_{k}),$$ where $D=G^{\mathfrak{S}}=G^{\mathfrak{U}}$.
  Let $U$ be a  submodular non-modular  subgroup of $G$ of minimal
order.   Then:}

(1) {\sl If $UU_{i}'/U_{i}'$ is modular in  $G/U_{i}'$ for
all $i=1, \ldots, k$, then $U$ is supersoluble.}

(2) {\sl If $U$ is  supersoluble and $UL/L$ is modular in  $G/L$ for
  all non-trivial nilpotent  normal subgroups $L$ of $G$, then
$U$ is a cyclic $p$-group for some prime $p$. }

{\bf Proof. }  Suppose that this lemma is false and let $G$ be a 
counterexample of minimal order.

(1)  Assume this is false. Suppose that 
$U\cap D\leq Z(D)$. Then every chief factor of $U$ below
 $U\cap Z(D)=U\cap D$ is cyclic  and, also,  $UD/D\simeq U/(U\cap  D)$ is
 supersoluble.
 Hence  $U$ is supersoluble, a contradiction. Therefore
 $U\cap D\nleq Z(D)$.   Moreover, Lemma 5.10(1)(2) 
 implies that $(U\cap D)Z(D)/Z(D)$ is
submodular in $D/Z(D)$ and so  $(U\cap
D)Z(D)/Z(D)$ is a non-trivial
normal  subgroup of  $D/Z(D)$  by Lemma 5.10(5).

  Hence for some $i$ we 
have $U_{i}/Z(D)\leq (U\cap 
D)Z(D)/Z(D),$ so  $U_{i}\leq (U\cap 
D)Z(D).$ But then $U_{i}'\leq  ((U\cap 
D)Z(D))'\leq U\cap D.$  By hypothesis, $UU_{i}'/U_{i}'=U/U_{i}'$ is modular  in  
$G/U_{i}'$ and so $U$ is modular in $G$ by \cite[p.~201, Property~(4)]{Schm}, a
 contradiction.
 Therefore Statement (1) holds.

(2) Assume that this is false.  Let $N=
U^{{\mathfrak{N}}}$ be the nilpotent residual of $U$.
 Then $N < U$ since $U$ supersoluble, so $N$ is modular in $G$.  
  It is clear that 
  every proper  subgroup $S$ of $U$ with 
$N\leq S$  is submodular in $G$, so the minimality of $U$ implies 
that $S$ is modular in $G$. Therefore, if $U$ has at least two distinct
  maximal subgroups $S$ and $W$
 such that $N\leq S\cap W$, then $U=\langle S, W \rangle $ is modular in 
$G$  by \cite[p. 201, Property (5)]{Schm}, contrary to our assumption on $U$.
 Hence $U/N$ 
is a cyclic $p$-group for some prime $p$ and $N\ne 1$ since $U$ is not cyclic.

Now we show that $U$ is  a $PT$-group. Let $S$ be a proper 
subnormal subgroup of $U$.  Then $S$ is  
submodular in $G$  since $U$ is  
submodular in $G$, so $S$ is modular  in $G$ and hence $S$ is 
quasinormal in $U$ by Theorem~3.1. Therefore 
 $U$ is a soluble $PT$-group, so $N=U^{{\mathfrak{N}}}=U'$ is a 
Hall abelian   subgroup of $U$ and every subgroup of $N$ is normal in $U$ by
 \cite[Theorem 2.1.11]{prod}.
Then $N\leq U^{{\frak{A}^{*}}}$ and so  $U^{{\frak{A}^{*}}}=NV,$ where  $V$ is a maximal
 subgroup of a Sylow $p$-subgroup $P\simeq U/N$ of $U$.  
Then 
$NV$ is modular in $G$ and $NV$ is subnormal in $G$ by  Lemma 5.10(4).  Therefore 
$NV$ is quasinormal in $G$ by Theorem~3.1.  Assume that for some minimal normal
 subgroup $R$ of $G$ we have  $R\leq (NV)_{G}$. Then $U/R$ is a modular in $G/R$
 by hypothesis, so  $U$ is modular in $G$, a contradiction. Therefore 
$(NV)_{G}=1$, so $NV$ is nilpotent  by \cite[Corollary 1.5.6]{prod} and
 then $V$ is normal in $U$.

Some maximal subgroup $W$ of $N$ is normal in $U$ with $|N:W|=q$. Then $S=WP$
 is a maximal subgroup of $U$ such that    $U/S_{U}$ is a non-abelian
 group of order $pq$.
  Hence $S$ is modular in $U$ by \cite[Lemma 5.1.2]{Schm}, so
 $S$ is modular in $G$. It follows that 
$U=NS$ is modular in $G$, a contradiction. 
 Therefore Statement  (2) holds.     

The lemma is proved.

{\bf Lemma 5.13.}   {\sl  If  $G$ is an
 $MT$-group, then  every quotient $G/N$ of $G$ is  an
$ MT$-group and  satisfies
 ${\bf M}_{P}$.}

{\bf Proof.}   
Let $L/N$ be submodular subgroup of $G/N$. Then $L$ is 
submodular subgroup in $G$ by Lemma 2.1(3), so $L$
 is modular in $G$ by hypothesis 
 and
 then $L/N$ is modular in $G/N$ by  \cite[p. 201, Property (3)]{Schm}. Hence 
$G/N$  is  an $MT$-group.
   
Now we show that $G/N$  satisfies
 ${\bf M}_{P}$.   Since $G/N$ is an  $ MT$-group, we can assume without loss of
 generality that $N=1$. Let $P/R$ be a normal non-abelian  $P$-subgroup 
  of $G/R$ and let $L/R\leq P/R$. Then $R/N$ is modular
 in $P/R$ by \cite[Lemma 2.4.1]{Schm}, so $L/R$ is submodular in $G/R$ and hence 
$L/R$ is modular in $G/R$ since $G/R$ is an $MT$-group. Therefore $G$  satisfies
 ${\bf M}_{P}$. 
 
The lemma is proved.

We need the following special case of Lemma 2.39.

{\bf Lemma 5.14.} {\sl  Suppose that $G$ has a Robinson
 complex $(D, Z(D);$ $ U_{1}, \ldots , U_{k})$ and let $N$ be a normal subgroup of $G$.  }

(1) {\sl If $N=U_{i}'$ and $k\ne 1$, then $Z(D/N)  =U_{i}/N =Z(D)N/N$ and 
 $$(D/N,
 Z(D/N); U_{1}N/N, \ldots , U_{i-1}N/N,
U_{i+1}N/N, \ldots
U_{k}N/N)$$ is  a Robinson complex of $G/N$.  }

(2) {\sl If $N$ is nilpotent, then  $Z(DN/N)= Z(D)N/N$ and
 $$(DN/N, Z(DN/N); U_{1}N/N, \ldots ,
U_{k}N/N)$$ is  a Robinson complex of $G/N$.}

 {\bf Proof of Theorem 5.8.}
 First assume that $G$ is an $M T$-group. Then $G$ is a $PT$-group and 
every quotient $G/N$ is 
an $M T$-group by Lemma  5.13. 
 Moreover, by Theorem 5.6,  
  $G$  has a normal perfect subgroup $D$ such that: 
  $G/D$ is a soluble group $PT$-group, and  
 if $D\ne 1$, $G$ has a Robinson complex
 $(D, Z(D); U_{1},  \ldots , U_{k})$ such that     
 for any set  $\{i_{1}, \ldots , i_{r}\}\subseteq \{1, \ldots , k\}$, where
 $1\leq r  < k$,  $G$ and $G /U_{i_{1}}'\cdots U_{i_{r}}'$ satisfy
 ${\bf N}_{p}$ for all $p\in \pi (Z(D))$ and
 ${\bf P}_{p}$  for all $p \in \pi (D)$. In view of  Lemma 5.13, $G$
 and $G /U_{i_{1}}'\cdots U_{i_{r}}'$ satisfy ${\bf M}_{p, q}$ for all pears
 $\{p, q\}\cap \pi (D)\ne \emptyset.$ 

In view of \cite[Theorem 2.1.11]{prod}, $G/D$ is a supersoluble $PT$-group
 and if $U/D$ is a 
 submodular subgroup of $G/D$, then $U$ is submodular in $G$ by Lemma 5.10(3), so
$U$ is modular in $G$ by hypothesis  and hence $U/D$ is modular in $G/D$ by 
\cite[p. 201, Property (4)]{Schm}, Therefore $G/D$ is an $M$-group by Proposition
 4.13.  
 Therefore the necessity
 of the condition of the theorem holds.

 Now, assume, arguing 
by contradiction, that $G$ is a non-$MT$-group of minimal order 
satisfying Conditions (i),  (ii), and (iii).  
 Then $D\ne 1$  
and   $G$ has a submodular subgroup $U$ such that $U$ is not modular
 in $G$ 
but  every submodular subgroup $U_{0}$   of $G$   with $U_{0} < U$   is modular in
$G$.   Let $Z=Z(D).$  Then $Z= \Phi (U_{i})=\Phi (D)$ since $D/Z$ is perfect.

(1) {\sl   $U$ is supersoluble.   }

 First assume that $k=1$, that is, $D=U_{1}=D'$.
 Then $(U\cap D)Z/Z$ is a submodular subgroup 
 of a simple group $D/Z$ by Lemma 5.10(1)(2). If $(U\cap D)Z/Z\ne 1$, then  
$(U\cap D)Z/Z=U_{i}/Z$  by Lemma 5.10(5), so  $(U\cap D)Z=U_{i}$, contrary to the fact 
that $Z\leq \Phi(U_{i})$.

 Hence $U\cap D\leq Z$,
 so $U\cap D=U\cap  Z$.  Therefore
 every chief factor of $U$  below $U\cap D$ is cyclic. On the other hand, $U/(U\cap D)
\simeq UD/D$ is supersoluble  by Condition (i)  since every $M$-group is
 supersoluble by  \cite[Theorem 2.4.4]{Schm}  and so $U$ is supersoluble.

Now let $k\ne 1$. We show that the hypothesis holds for $G/U_{i}'$
 for all $i=1, \ldots , k$. We can assume   without loss of generality that $i=1$.
  Let $N=U_{1}'$.  Then
  $ (G/N)/(D/N)\simeq G/D$ is an $M$-group   and $
(D/N)'=D '/N=D/N$. From Lemma
 5.14(1) it follows that   
 $(D/N, Z(D/N); U_{2}N/N, \ldots , U_{k}N/N)$ is
 a Robinson complex of $G/N$ and 
  $$Z(D/N)=U_{1}/N=ZN/N\simeq Z/(Z\cap N).$$

   Moreover, if
 $\{i_{1}, \ldots , i_{r}\}\subseteq \{2, \ldots , k\}$, where $2\leq r  < k$,
 then the quotients
 $G/N=G/U_{1}'$ and $$(G/N) /(U_{i_{1}}N/N)'\cdots (U_{i_{r}}N/N)'=
(G/N)/(U_{i_{1}}'\cdots U_{i_{r}}'U_{1}'/N)\simeq G/U_{i_{1}}'\cdots U_{i_{r}}'U_{1}'$$
 satisfy
 ${\bf N}_{p}$ for all
 $p\in  \pi (Z(D/N))\subseteq \pi (Z)$, ${\bf P}_{p}$ for all
 $p\in  \pi (D/N))\subseteq \pi (D)$, and  ${\bf M}_{p, q}$  for all pears
  $\{p, q\}\cap \pi (D)\ne \emptyset$  by Condition (iii).  

Therefore the hypothesis holds for $G/N=G/U_{1}'$, so the 
submodular subgroup $UU_{1}'/U_{1}'$ of $G/U_{1}'
$  is modular in
  $G/U_{1}'$ by the choice of $G$.  Hence $U$ is  supersoluble by
Lemma 5.12(1).  
  
(2) {\sl   If  $N\ne 1$
 is a normal nilpotent subgroup of $G$, then Conditions (i), (ii), and (iii) 
hold for $G/N$ } (See the proof of Claim (1) and use the fact that every quotien of an
 $M$-group is an $M$-group as well).

(3) {\sl If $ XN/N$  is submodular in $G/N$ for some
 normal nilpotent subgroup $N\ne 1$, then $ XN/N$  is modular in $G/N$ and 
 $XN$ is modular in $G$. In particular, $U_{G}=1$ }  (This follows from
 Claim (2), Lemma 5.10(3) and the choice of $G$).

 (4) {\sl $U$ is a cyclic $p$-group for some prime $p$ and 
 $U\cap  Z_{\infty}(G)$ is the maximal subgroup of $U$.}

Let $N$  be  a nilpotent 
 non-identity  normal subgroup of $G$. Then $UN/N$ is submodular in $G/N$
 by Lemma 5.10(2), so $UN/N$ is modular in $G/N$ by Claim (3). Hence $U$ is a cyclic
 $p$-group for some prime $p$ by Lemma 5.12(2) and Claim (1).

Now let $V$ be a maximal subgroup of $U$. Then $V=U^{{\frak{A}^{*}}}$  
  is subnormal in $G$ by Lemma 5.10(4), hence  $V$ is  modular  and so
quasinormal  in $G$ by Theorem A. Therefore 
 $V\leq Z_{\infty}(G)$ by \cite[Corollary 1.5.6]{prod} since $V_{G}=1=U_{G}$
 by Claim (3). Hence $V=U\cap  Z_{\infty}(G)$.

(5) {\sl $G$ has a normal subgroup $C_{q}$ of order  $q$ for some  $q\in \pi (Z(D))$. }

It is enough to show that $Z\ne 1$. Assume
 that this is false. Then  $D= U_{1}\times \cdots 
\times U_{k}$, where $U_{i}$ is a simple non-abelian group for all $i$.  

Let $E=U_{i}U$. 
 We show that $U_{i}\leq N_{G}(U) $. Note that $U$ is submodular in $E$
 by Lemma 5.10(1). Therefore there is a subgroup chain 
$$U=E_{0} < E_{1} < \cdots < E_{t-1} < E_{t}=E$$ such that $E_{i-1}$  is a maximal
modular subgroup of  $E_{i}$ for all $i=1, \ldots, t$ and for $M=E_{t-1}$ we have 
$M=U(M\cap U_{i})$. Moreover,  by  \cite[Lemma 5.1.2]{Schm}, either
 $M=E_{t-1}$ is  a maximal normal subgroup of $E$ or $M$
 is a maximal subgroup of $E$ such  that 
$E/M_{E}$ is a non-abelian group of order $qr$ for primes  $q$ and $r$.
In the former case we have $M\cap U_{i}=1$, so $U=M$ is normal in $E$. The second case
 is impossible 
since $E$ has no a quotient of order $qr$. Therefore  $U_{i}\leq N_{E}(U)$ for all $i$,
 so 
$D\leq N_{G}(U)$ and hence $U\cap D\leq O_{p}(D)=1$. It follows than $DU=D\times U$, so 
$1  < U\leq C_{G}(D)$. But  $C_{G}(D)\cap D=1$ since $Z(D)=1$. Therefore
 $C_{G}(D)\simeq 
C_{G}(D)D/D$ is soluble. Hence for some prime $q$ dividing $|C_{G}(D)|$ we have 
$O_{q}(C_{G}(D))\ne 1$. But $O_{q}(C_{G}(D))$ is characteristic in $C_{G}(D)$, so 
 we have (5).  

(6) {\sl $U^{G}$ is soluble}.

The subgroup $C_{q}U/C_{q}$  is submodular in $G/C_{q}$ by Lemma 5.10(2), so 
$C_{q}U$  is modular in $G$ by Claim (3).

 First assume that $C_{q}U$  is not subnormal  in $G$.
 Then 
$C_{q}U$  is not quasinormal in $G$, so 
$(C_{q}U)^{G}/(C_{q}U)_{G}$  is a non-abelian   $P$-group  by Lemma 5.11. 
Hence $(C_{q}U)^{G}$  is soluble.

Now assume that  $C_{q}U$  is  subnormal  in $G$,  so 
$C_{q}U/C_{q}$  is a subnormal $p$-subgroup of   $G/C_{q}$ and hence 
 $$U^{G}/(U^{G}\cap C_{q})\simeq
 C_{q}U^{G}/C_{q}=(C_{q}U/C_{q})^{G/C_{q}}\leq O_{p}(G/C_{q}).$$
 Hence  $U^{G}$ is soluble.

(7) {\sl $U$ is not subnormal in $G$. }

Assume that   $U$ is subnormal in $G$.  Then $U$ is quasinormal in $G$ since 
$G$ is a $PT$-group by Theorem 5.6, so $U$ is modular in $G$ by Theorem 3.1,
 a contradiction. Hence we have (7).

(8)  $|U|=p$.

Assume that  $|U| > p$. Then $1 < V\leq R:=O_{p}(Z_{\infty}(G))$ by Claim (4) and 
  $U\nleq R$ by Claim (7). Let $E=RU$. Then $E$ is not subnormal in $G$ by Claims 
(4) and (7) and this subgroup is modular in $G$ by Claim (3) and Lemma 5.10(2). Moreover,
   $UR/R\simeq U/(U\cap R)=U/V$  has order $p$, so $(RU)_{G}=R$.
    Therefore,
 in view of Lemma 5.11, $$G/R=E^{G}/E_{G}\times K/E_{G}=U^{G}R/R\times K/R,$$  where 
 $RU^{G}/R\simeq U^{G}/(U^{G}\cap R)$ 
 is a   non-abelian  $P$-group of order prime to $K/R$. Then $RU^{G}/R$ is a $\pi$-group,
where $\pi =\{p, q\}$  for some prime $q$, so $G$ is $\pi$-soluble and hence 
  $D$ and $D/Z$ are $\pi$-soluble groups.

 Assume that  $U^{G}\cap D\ne 1$. Since 
 $U^{G}\cap D\leq Z=Z(D)\leq \Phi (D)$ by Claim (6),
 for some $i$ and for some  $r\in \{p, q\}$ the mumber $r$ divides $|U_{i}/Z|$.
 It follows that $U_{i}/Z$ is an abelian group, a contradiction.    
Therefore $U^{G}\cap D= 1$, so  from $$U^{G}\simeq U^{G}/(U^{G}\cap D)\simeq U^{G}D/D$$
we get that   $U^{G}$ is an $M$-group. 
 Then from Lemma 2.4.1 and Theorem 2.4.4 in \cite{Schm} it follows that  
$ U(U^{G}\cap R)/(U^{G}\cap R)$ is a Sylow $p$-subgroup of 
$U^{G}/(U^{G}\cap R)$ and  so  $U(U^{G}\cap R)$ is a cyclic Sylow $p$-subgroup of $U^{G}$.
It follows that either  $U(U^{G}\cap R)=U$ or $U(U^{G}\cap R)=U^{G}\cap R$. In the former case
 we have  $U^{G}\cap R=V$, which is impossible by  Claim (3), so 
$U(U^{G}\cap R)=U^{G}\cap R$ and hence $U$ is subnormal in $G$, contrary to Claim (7).
Therefore  we have (8).

(9) {\sl $U\nleq D$.}

Assume $U\leq D$. 
From Claim (7) it follows that $U\nleq Z$ and then, by Claim (8) and
 Lemma 5.10(1)(2)(5), for some $i$ we have 
 $U\simeq UZ/Z=U_{i}/Z$, so $UZ=U_{i}$, a
 contradiction.  Hence we have~(9).

(10) $O_{p}(D)=1$.

Assume that $G$ has a normal subgroup $Z_{p}\leq Z$ of order $p$.
Then $Z_{p}U$ is not subnormal in $G$ by Claim (7) and, also,  $(Z_{p}U)_{G}=Z_{p}$
 by Claim (8) and 
$(Z_{p}U)^{G}=Z_{p}U^{G}$, so 
$$G/Z_{p}=Z_{p}U^{G}/Z_{p}\times K/Z_{p},$$ where $Z_{p}U^{G}/Z_{p}$ is
 a non-abelian $P$-group of order $p^{a}q^{b}$ prime to $|K/Z_{p}|$. Hence 
$G/Z_{p}$,   $D/Z_{p}$, and $D$  are $\{p, q\}$-soluble, where $p$ divides $|D/Z_{p}|$. 
Hence $O_{p}(D/Z)\ne 1$. This contradiction completes the proof of the  claim.

(11) {\sl $DU=D\times U$.   In particulart, $NU$ is not subnormal in $G$ for every normal
 subgroup $N$ of $G$ contained in $D$.}

In view of Claims (8) and (9), it is enough to show that $U_{i}\leq N_{G}(U)$ for all $i$.

Let $E=U_{i}U$ and let 
$$U=E_{0} < E_{1} < \cdots < E_{t-1} < E_{t}=E$$ be  a   subgroup chain such that $E_{i-1}$
  is a maximal
modular subgroup of  $E_{i}$ for all $i=1, \ldots, t$ and for $M=E_{t-1}$ we have 
$M=U(M\cap U_{i})$ and either
 $M=E_{t-1}$ is  a maximal normal subgroup of $E$ or $M$
 is a maximal subgroup of $E$ such  that 
$E/M_{E}$ is a non-abelian group of order $qr$ for some primes  $q$ and $r$.

First assume that $M$ is normal in $E$.
 From $E=U_{i}U=U_{i}M$ it follows that $E/M\simeq U_{i}/(M\cap U_{i})$ is a
 simple group, so  $M\cap U_{i}=Z$ and hence $U\cap U_{i}=1$ by Claim (7).
 Then, by the Frattini argument,
 $E=MN_{E}(U)=ZN_{E}(U)$.
 But $Z\leq \Phi (E)$ since $Z\leq \Phi (U_{i})$. Therefore  $N_{E}(U)=E$,
 so $U_{i}\leq N_{G}(U)$.

Finally, assume that $E/M_{E}$ is a non-abelian group of order $qr$. Then
 $U_{i}/(U_{i}\cap M_{E})$ is  soluble, so $U_{i}=U_{i}\cap M_{E}$ 
 since $U_{i}$ is perfect. Therefore $U_{i}, U\leq M$ and so $M=E $, 
 a contradiction. Hence we have (11).

(12) {\sl $U^{G}$ is not a non-abelian $P$-group. }

Assume that   $U^{G}=Q\rtimes U$ is a non-abelian $P$-group of type $(q, p)$
 and let $\pi =
\{q, p\}$.  Let $ S=U^{G}\cap D  $. In view of Claim (9) and Lemma 
 2.2.2 in  \cite{Schm},  $U\nleq S$ and so  $S\leq O_{q}(D)\leq Z(D)$.
 Moreover, Claim (11) implies that 
$S\ne Q$, so  $DU^{G}/D\simeq U^{G}/S$ is a 
 non-abelian $P$-group of type $(q, p)$ by \cite[Lemma 2.2.2]{Schm}}. In view of
 Claims (8) and  (11), $(DU)_{G}=D$  and hence  
$G/D=DU^{G}/D \times  K/D,$ where 
$DU^{G}/D=O_{\pi}(G/D)$ and $K/D=O_{\pi'}(G/D)$, so $G/D$ is $\pi$-decomposable.

Let $E$ be a minimal supplement to $N$ in $G$. Then $E\cap D\leq \Phi(E)$, so $E$ is
 soluble and $\pi$-decomposable, that is, $E=O_{\pi}(E)\times O_{\pi'}(E)$ since
 $G/D\simeq E/(E\cap D)$.

 First suppose that  $\pi \cap \pi (D)= \emptyset$ and 
    let $G_{r}$ be a Sylow $r$-subgroup of
 $E$, where $r\in \pi$.  Then  $G_{r}$ is a Sylow $r$-subgroup of
 $G$, therefore $E$ has a Hall $\pi$-subgroup $H$ since  $E$ is soluble
 and $U^{G}=Q\times U\leq H^{x}$ for all $x\in G$. 
 Then,  by Condition (i), 
for every $r$-element $x$ of $G$, where $r\in \pi$, $U$ is modular in 
$\langle x, U \rangle$ since $H\simeq DH/D$ is an $M$-group.

Now let $x\in G_{r}$, where  $r\not \in  \pi$. Then for some Sylow
 $r$-subgroup $D_{r}$ of $D$ 
and a Sylow $r$-subgroup $E_{r}$ of $E$ and some $y\in G$ we have 
 $G_{r}=D_{r}E_{r}^{y}$. Hence $x=de$, where  $d\in D_{r}$ and $e\in E_{r}^{y}$. Now note 
that for any $u\in U$ we have $d^{-1}u^{-1}du\in D\cap U^{G}=1$, so $D\in C_{G}(U)$.
On the other hand, $e$ is a $\pi'$-element of the  $\pi$-decomposable group $E^{y}$, so 
$e\in C_{G}(U)$. Therefore  $x\in C_{G}(U)$ and hence $U$ is normal in 
 $\langle x, U \rangle$.  Therefore $U$ is modular in $G$ by 
  \cite[Theorem 5.1.13]{Schm},
 a contradiction. Finally, if 
  $\pi \cap \pi (D)\ne  \emptyset$, then $U$ is modular in $G$ by Condition (iii).
This contradiction completes the proof of the claim.

{\sl Final contradiction.} From Claims (5), (7), (9), and (11) it follows that
 $E=C_{q}U=C_{q}\times U$ is not subnormal in $G$ and, in view of Claim (8),
 $E_{G}=C_{q}$.
Hence $G/E_{G}= E^{G}/E_{G}\times K/E_{G},$ where
 $$E^{G}/E_{G}=C_{q}U^{G}/C_{q}\simeq 
U^{G}/(C_{q}\cap U^{G})$$ is a non-abelian
 $P$-group of order prime to $|K/C_{q}|$ by Lemma 5.11. Hence $G$ is a $\pi$-soluble 
group, where $\pi= \pi (U^{G}/(C_{q}\cap U^{G}))$. Then $D/C_{q}$ is 
$\pi$-soluble. But $C_{q}\leq \Phi (D)$, so $q$ divides $|D/C_{q}|$. Hence
$q$ does not divides  $|C_{q}U^{G}/C_{q}|$.

If $C_{q}\cap U^{G}=1$, then $U^{G}\simeq  C_{q}U^{G}/C_{q} $ is a non-abelian
 $P$-group, contrary to Claim (12),  so  $C_{q}\leq  U^{G}$. Then   
 $C_{q}$ is a Sylow $q$-subgroup of $U^{G}$.
 Hence $U^{G}=C_{q}\rtimes (R\rtimes U)$,
 where  $R\rtimes U\simeq  U^{G}/C_{q}$ is a non-abelian $P$-group. 
Let $C=C_{U^{G}}(C_{q})$. Then $U\leq C$ by Claim (11) and so, by Lemma 
 2.2.2 in  \cite{Schm},  $R\rtimes
 U=U^{R\rtimes U}\leq C$. Hence $C_{q}\leq Z(U^{G})$.
 Therefore
  $U^{G}=C_{q}\times (R\rtimes U)$, where $R\rtimes U$ is characterisric in $U^{G}$
 and so it is normal in $G$. But then $U^{G}=R\rtimes U\ne C_{q}\rtimes (R\rtimes U)$, 
a contradiction. 
              
The theorem is proved.

   \end{document}